\documentclass[aap,nameyear,MSNbibl,seceqn,dvips]{arximspdf}

\doi{10.1214/10-AAP747}
\volume{21}
\issue{5}
\pubyear{2011}
\firstpage{1965}
\lastpage{1993}

\makeatletter
\newtheorem{theorem}{Theorem}[section]
\newtheorem{propos}[theorem]{Proposition}
\newproclaim{Remark}[theorem]{Remark}
\newproclaim{rem}[theorem]{Remark}

\newcommand{\1}{1}
\newcommand{\N}{\mathbb{N}} 
\newcommand{\R}{\mathbb{R}} 
\newcommand{\esssup}{\operatorname{ess\sup}}
\newcommand{\autorF}{\citet{Faller-09}}
\newcommand{\autorFx}[1]{\citet{Faller-09}, #1}
\newcommand{\autorFR}{\citet{FallerRueschendorf-09}}
\newcommand{\autorFRx}[1]{\citet{FallerRueschendorf-09}, #1}
\newcommand{\autorKR}[1]{\citeauthor{KuehneRueschendorf#1}}
\newcommand{\autorKRx}[1]{\citet{KuehneRueschendorf#1}}
\newcommand{\autorKRxx}[2]{\citet{KuehneRueschendorf#1}, #2}

\newcommand{\refeq}[1]{\textup{(\ref{#1})}} 
\renewcommand{\epsilon}{\varepsilon}
\newcommand{\eqref}[1]{\textup{(\ref{#1})}}

\newcommand{\sideset}[3]{{}#1#3#2}
\renewcommand{\mid}{|}
\makeatother

\begin{document}
\begin{frontmatter}

\title{On approximative solutions of multistopping~problems}
\runtitle{On approximative solutions of multistopping problems}

\begin{aug}
\author[A]{\fnms{Andreas} \snm{Faller}\ead[label=e1]{afaller@hotmail.com}}
\and
\author[A]{\fnms{Ludger} \snm{R\"{u}schendorf}\corref{}\ead[label=e2]{ruschen@stochastik.uni-freiburg.de}\ead
[label=u1,url]{http://www.stochastik.uni-freiburg.de}}
\runauthor{A. Faller and L. R\"{u}schendorf}
\affiliation{University of Freiburg}
\address[A]{
Department of Mathematical Stochastics\\
University of Freiburg\\
Eckerstr. 1\\
79104 Freiburg\\
Germany\\
\printead{e1}\\
\hphantom{\textsc{E-mail}: }\printead*{e2}\\
\printead{u1}} 
\end{aug}

\received{\smonth{3} \syear{2010}}
\revised{\smonth{11} \syear{2010}}

%
\begin{abstract}
In this paper, we consider multistopping problems for finite discrete
time sequences $X_1, \ldots , X_n$.
$m$-stops are allowed and the aim is to maximize the expected value of
the best of these $m$ stops.
The random variables are neither assumed to be independent not to be
identically distributed. The basic assumption
is convergence of a related imbedded point process to a continuous time
Poisson process in the plane,
which serves as a limiting model for the stopping problem. The optimal
$m$-stopping curves for this limiting model are
determined by differential equations of first order. A general
approximation result is established which ensures
convergence of the finite discrete time $m$-stopping problem to that in
the limit model. This allows the construction of approximative
solutions of the discrete time $m$-stopping problem. In detail, the
case of i.i.d. sequences with
discount and observation costs is discussed and explicit results are obtained.
\end{abstract}

%

\end{frontmatter}

\section{Introduction}\label{sec:1}

In this paper, we consider multistopping problems for discrete time
sequences $X_1,\ldots ,X_n$. In comparison to the usual stopping
problem, there are $m$ stops $1 \le T_1 <\cdots <T_m \le n$ allowed. The
aim is to determine these stopping times in such a way that
%
\begin{equation}\label{eq:1.1}
E\Bigl[ \max_{1\le i\le m} X_{T_i}\Bigr] = E[X_{T_1}\vee\cdots \vee X_{T_m}] =
\mathrm{sup}.
\end{equation}
Thus, the gain of a stopping sequence $(T_i)_{i\le m}$ is the expected
maximal value of the $m$ choices $X_{T_i}$.
In the case $m=1$, this stopping problem reduces to the classical Moser
problem [\citet{Moser-56}]. We will see that optimal
$m$-stopping times exist and are determined by a recursive description.

Our aim is to obtain \textit{explicit} approximative solutions of the
$m$-stopping problem in (\ref{eq:1.1}) under general distributional
conditions. In particular, we do not assume that the random variables
$X_i$ are independent or identically distributed or
are even of specific i.i.d. form with $X_i \sim U(0,1)$ as assumed in
several papers in the literature. Our basic assumption is convergence
of the imbedded planar point process (\ref{eq:1.2}) of rescaled
observations to some Poisson point process $N$
in the plane,
%
\begin{equation}\label{eq:1.2}
N_n = \sum_{i=1}^{n} \delta_{( {i}/{n}, X^{n}_i)} \stackrel{d}{\rightarrow} N.
\end{equation}
Here $X^n_i = \frac{X_i - b_n}{a_n}$ is a nomalization of the $X_i$
induced typically from the central limit theorem for
maxima respectively related point process convergence results.
Our aim is to prove that under some regularity conditions the optimal
$m$-stopping problem of $X_1, \ldots , X_n$
can be approximated by a suitable formulated $m$-stopping problem for
the continuous time Poisson process $N$ which
serves as a limiting~model for the discrete time model. Furthermore, we
want to show that the stopping problem in the limit
model can be solved in explicit form. The solution is described by an
increasing sequence of stopping curves with their related
threshold stopping times. These curves solve usual one-stopping
problems for transformed Poisson processes and are characterized by
differential equations of first order, which can be solved either in
exact form or numerically.
The solution for the limit model also allows us to construct approximative
optimal stopping times
for the discrete time model. We apply this approach in detail to the
$m$-stopping of sequences
$X_i= c_i Z_i + d_i$
\textit{with discount} and observation costs and i.i.d. sequences $Z_i$.

It has been observed in several papers in the literature that optimal
stopping may have an easier
solution in a related form for a Poisson number of points or for
imbedded homogeneous Poisson processes
as for instance in the classical house selling problem or in best
choice problems. For $m=1$
[see, e.g., \citet{ChowRobbinsSiegmund-71}, \citet{Sakaguchi-76},
\citet{BrussRogers-91}, \citet{GnedinSakaguchi-92},
\citet{Gnedin-96}, \citet{BaryshnikovGnedin-00}].
For general reference, we refer to \citet{Ferguson-07}, Chapter 2.
For $m \geq1 $, multistopping problems were introduced in \citet
{Haggstrom-67} who derived some
structural results corresponding roughly to Theorem \ref{theo:2.2};
compare also some extensions in \citet{Nikolaev-99}.
The two stopping problem has been considered in the case of Poissonian
streams in \citet{SaarioSakaguchi-92}. In
this paper, differential equations were derived corresponding to those
for the one-stopping problems in \citet{Karlin-62},
\citet{Siegmund-67} and \citet{Sakaguchi-76}. Multiple buying---selling problems were studied in \citet{BrussFerguson-97}
based on a vector valued formulation with pay-off given by the sum of
the $m$-choices instead of the max as in (\ref{eq:1.1});
see also the extension in \citet{Bruss-10}.
In \autorKR{-02}\ (\citeyear{KuehneRueschendorf-02}) the case of 2-stopping problems for i.i.d. sequences
was treated based on the approximative approach
in \autorKR{-00a}\ (\citeyear{KuehneRueschendorf-00a}).
The results in this paper were rederived in \citeauthor
{AssafGoldsteinCahn-04} (\citeyear{AssafGoldsteinCahn-04,AssafGoldsteinCahn-06}) and in \citet{GoldsteinCahn-06}.
In case $m=1$ based on this approximation for several classes of
independent and dependent sequences
optimal solutions have been found in explicit form [see \autorKR{-00b}\
(\citeyear{KuehneRueschendorf-00b,KuehneRueschendorf-04}) and \autorFR].
The present paper establishes an extension of the approximative
approach as described above to $m$-stopping problems as in
(\ref{eq:1.1}). It is based on the dissertation of \citet{Faller-09}
to which we refer for some technical details in the proofs.

The program to establish this approximation approach in general is
based on the following steps.
In Section \ref{sec:2}, we formulate the necessing recursive
characterization of the optimal solutions of the $m$-stopping problem
corresponding to Bellman's optimality equation. Section~\ref{sec:3} is devoted to
solve the $m$-stopping problem for the
limit model of an inhomogeneous Poisson process. A particular
difficulty arises from the fact that in the limit model
the intensity function is typically infinite along a lower boundary
curve, In consequence, known stationary Markovian
techniques as for homogenous Poisson processes do not apply. The main
result, Theorem~\ref{theo:3.3}, shows that
the optimal $m$-stopping problem can be reduced to $m$ 1-stopping
problems for transformed Poission processes.
The optimal stopping curves are characterized by a sequence of
differential equations of first order.

In Section \ref{sec:4}, we are able to derive explicit solutions for
some classes of differential equations, as appearing in the description of
the optimal stopping curves. This part is based on developments in
\autorFR\ for the
case $m=1$. Section~\ref{sec:5} gives the basic approximation theorem
(Theorem~\ref{theo:5.2}) allowing to
approximate the finite discrete problems by $m$-stopping in the limit
model. The proof of this result needs to develop a new technique. It
also uses essentially the extension of the convergence of multiple
stopping times in Proposition~\ref{prop:5.1}
in \autorFR\ for $m=1$ to $m \geq1$.
We restrict our presentation to the essential new part of this proof.
Finally in Section~\ref{sec:6} we obtain as
application solutions in \textit{explicit} form for optimal $m$-stopping
problems for sequences $X_i= c_i Z_i + d_i$
with $Z_i$ i.i.d. and with discount and observation costs $c_i$, $d_i$. It
is remarkable that we get detailed results including
the asymptotic constants as well as approximative optimal stopping
sequences in explicit form. Our aim is to extend these results in
subsequent papers to further classes of stopping problems as to
selection problems, to the sum cost case as well
as to some classes of dependent sequences. It seems also possible as
done in the case $m=1$, to extend this approach to
the case where cluster processes arise in the limit.

\section{$m$-stopping problems for finite sequences}\label{sec:2}

In this section, we give a formulation of the optimality principle for
the $m$-stopping of discrete recursive sequences.
Given a discrete time sequence $(X_i,\mathcal F_i)_{1\le i\le n}$ in a
probability space $(\Omega,\mathcal A,P)$ with filtration $\mathcal F
= (\mathcal F_i)_{0\le i\le n}$ the $m$-stopping problem ($1\le m\le
n$) is to find stopping times $1\le T_1< T_2<\cdots<T_m\le n$ w.r.t.
the filtration $(\mathcal F_i)_{1\le i\le n}$ such that
%
\begin{equation}\label{eq:2.1}
E\Bigl[ \max_{1\le i\le m} X_{T_i} \Bigr]= E[ X_{T_1}\vee\cdots \vee X_{T_m}] =
\mathrm{sup}.
\end{equation}

In case $m=1$, \eqref{eq:2.1} is identical to the usual (one-)stopping
problem. A well-known recursive solution of this problem [see
\citet{ChowRobbinsSiegmund-71},
 Theorem 3.2] is based on the threshold curves
$W_i=W_F(X_{i+1},\ldots ,X_n)$ of the optimal stopping time defined by
%
\begin{eqnarray}\label{eq:2.2}
W_n &:=& -\infty, \nonumber
\\[-8pt]
\\[-8pt]
W_i &:=& E [X_{i+1}\vee W_{i+1} \mid \mathcal{F}_{i}] \qquad\mbox{for } i=n-1,\ldots ,0.
\nonumber
\end{eqnarray}
We need a version of this classical result for stopping times larger
than a given stopping time $S$.
%

%
\begin{propos}[(Recursive solution of one-stopping problems)]
\label{prop:2.1}
{\smallskipamount=0pt
\begin{longlist}[(b)]
\item[(a)] For any time point $0\le k\le n-1$, the $\mathcal
F$-stopping time
\[
T(k):=\min\{k<i\le n \dvtx  X_i>W_i\}
\]
is optimal in the sense that for any $\mathcal F$-stopping time $T>k$
we have
%
\begin{equation}\label{eq:2.3}
E\bigl[X_{T(k)}\mid\mathcal F_k\bigr] = W_k\ge E[X_T\mid\mathcal F_k] \qquad P
\mbox{-a.s.}
\end{equation}

\item[(b)]
For any $\mathcal F$-stopping time $S$, the $\mathcal F$-stopping time
\[
T(S)=\min\{S<i\le n \dvtx  X_i>W_i\}
\]
is optimal in the sense that for any $\mathcal F$-stopping time $T$
with $S<T$ on $\{S<n\}$ and $S=T$ on $\{S=n\}$ we have
%
\begin{equation}
E\bigl[X_{T(S)}\mid\mathcal{F}_{S}\bigr] = W_S\ge E[ X_T \mid\mathcal{F}_{S}]
\qquad P \mbox{-a.s.}
\end{equation}
\end{longlist}}
\end{propos}

\begin{Remark}\label{rem:2.1}
For $m$ stopping problems, the following variant of Proposition~\ref
{prop:2.1} will also be needed
[for details of the proof, see \autorF].

Let $Y_1,\ldots , Y_n\dvtx (\Omega,\mathcal A,P)\to E$ be random variables
taking values in a measurable space $E$ and $\mathcal F:=(\mathcal
F_i)_{0\le i\le n}$ a filtration in $\mathcal A$ such that $\sigma
(Y_i)\subset\mathcal F_i$ for all $1\le i\le n$. Let $S$ be an
$\mathcal F$-stopping time, let $Z\dvtx (\Omega,\mathcal A,P)\to\overline
\R$ be $\mathcal F_S$-measurable and $h\dvtx E\times\overline\R\to
\overline\R$ be measurable with $Eh(Y_i,Z)^+<\infty$. Also define
recursively for $z\in\overline\R$
%
\begin{eqnarray} \label{eq:2.5a}
W_n(z) &:=& h(Y_n,z),\nonumber
\\[-8pt]
\\[-8pt]
W_i(z) &:=& E[h(Y_{i+1},z)\vee W_{i+1}(z) \mid \mathcal{F}_{i}]
\qquad\mbox{for } i=n-1,\ldots ,0.
\nonumber
\end{eqnarray}
Then the $\mathcal{F}$ stopping time
%
\begin{equation} \label{eq:2.6a}
T(S,Z) := \min\{S<i\le n \dvtx  h(Y_i,Z) > W_i(Z_i) \},
\end{equation}
where $Z_i:=Z \1_{\{S\le i\}}$ is optimal in the sense that for any
further $\mathcal{F}$-stopping time $T$ with $S<T$ on $\{S<n\}$ and
$S=T$ on $\{S=n\}$ we have
%
\begin{equation} \label{eq:2.7a}
E \bigl[h\bigl(Y_{T(S,Z)},Z\bigr)|\mathcal{F}_{S} \bigr] %
= W_S(Z_S) \ge E  [ h(Y_T,Z)|\mathcal{F}_{S}  ] \qquad P \mbox{-a.s.}
\end{equation}
\end{Remark}

Similar as for the one-stopping problems the idea of solving \eqref
{eq:2.1} is simple. The $\ell$th stopping time $T_\ell$ should be $i$
if the $(m-\ell)$-stopping value past $i$ with guarantee value $X_i$
is in expectation larger than the $(m-\ell+1)$-stopping value past $i$
and with guarantee value reached before time $i$. %
This idea leads to the following construction. Define $W_i^0(x):=x$ for
$x\in\overline\R$ and inductively for $1\le m\le n$, $x\in\R$
define thresholds $W_k^m(x)$ by
%
\begin{eqnarray}\label{eq:2.5}
 \qquad W^m_{n-m+1}(x)
& :=& x,\nonumber
\\[-8pt]
\\[-8pt]
 \qquad W^m_i(x)
&:=& E[ W^{m-1}_{i+1} (X_{i+1}) \vee W^m_{i+1}(x) \mid \mathcal
{F}_{i}] \qquad\mbox{for } i=n-m,\ldots ,0.
\nonumber
\end{eqnarray}
The related threshold stopping times are defined recursively for $k\le
n-m$ by
%
\begin{eqnarray} \label{eq:2.6}
T^{m}_{1} (k,x) &:=& \min\{k< i\le n-m+1 \dvtx  W^{m-1}_i (X_i)>W^m_i(x)\},\hspace*{-25pt}
\nonumber
\\
T^{m}_{\ell}(k,x) &:=& \min\{T^{m}_{\ell-1}(k,x)<i\le n-m+\ell\dvtx\hspace*{-25pt}\\
&&\hphantom{\min\{}    W^{m-l}_i(X_i)> W^{m-l+1}_i(x\vee M_{\ell
-1, i}) \}\hspace*{-25pt}
\nonumber
\end{eqnarray}
for $2\le\ell\le m$ and $M_{j,i} := X_{T_j^m (k,x)} \1_{\{
T_j^m(k,x)\le i\}}$.

Equation \eqref{eq:2.6} corresponds to a sequence of $m$ one-stopping problems
for (more complicated) transformed sequences of random variables. The
following result extends the classical recursive characterization of
optimal stopping times for one-stopping problems in Proposition \ref
{prop:2.1} to the case $m\ge1$. Related structural results can be
found in the papers of \citet{Haggstrom-67}, %
\citet{SaarioSakaguchi-92}, %
\citet{BrussFerguson-97}, %
\citet{Nikolaev-99}, \citet{BrussDelbaen-01} and \autorKR{-02}\ (\citeyear
{KuehneRueschendorf-02}).

\begin{theorem}[(Recursive characterization   of   $m$-stopping problems)]%
\label{theo:2.2}
The $\mathcal F$-stopping times $(T_\ell^m(k,x))_{1\le\ell\le m}$
are optimal in the sense that for all $\mathcal F$-stopping times
$(T_\ell)_{1\le\ell\le m}$ with $k<T_1<\cdots <T_m\le n$ we have
\begin{eqnarray*}
 &&E\bigl[ x\vee X_{T^{m}_1(k,x)}\vee\cdots \vee X_{T^{m}_m(k,x)}
\mid\mathcal{F}_k\bigr]   \\
&& \qquad = E \bigl[ W^{m-1}_{T^{m}_{1}(k,x)} \bigl(x\vee X_{T^{m}_{1}(k,x)}\bigr) \mid
\mathcal{F}_k\bigr] = W^m_k(x)
\\
&& \qquad \ge E [ x\vee X_{T_1}\vee\cdots \vee X_{T_m} \mid\mathcal{F}_k]
\qquad P \mbox{-a.s.}
\end{eqnarray*}
\end{theorem}

The proof of Theorem \ref{theo:2.2} follows by induction in $m$
based on Proposition \ref{prop:2.1} and Remark \ref{rem:2.1}
similarly as in the case $m=1$. For details, see \autorFx{Satz~2.1}\
or \autorKRxx{-02}{Proposition 2.1}. In general, the recursive
characterization of optimal $m$-stopping times and values is difficult
to evaluate. Our aim is to prove that one can construct optimal
$m$-stopping times and values approximatively by considering related
limiting $m$-stopping problems for Poisson processes in continuous time.

\section{$\mathbf m$-stopping of Poisson processes}
\label{sec:3}

In this section, we deal with the optimal $m$-stopping problem for the
limit model given by a Poisson point process $N$.
We consider a Poisson process $N=\sum_k\delta_{(\tau_k,Y_k)}$ in the
plane restricted to some set
\[
M_f=\{(t,x)\in[0,1]\times\overline\R; x>f(t)\},
\]
where $f\dvtx [0,1]\to\R\cup\{-\infty\}$ is a continuous lower boundary
function of $N$. The intensity of $N$ may be
(and in typical cases is)
infinite along the lower boundary~$f$. As in \autorKRx{-00a},
respectively,
\autorFR\ who consider the case $m=1$, we assume that the intensity
measure $\mu$ of $N$ is a Radon measure on $M_f$ with the topology on
$M_f$ induced by the usual topology on $[0,1]\times\overline\R$.
Thus any compact set $A\subset M_f$\vspace*{-2pt} has only finitely many points. By
convergence in distribution ``$N_n \stackrel{d}{\to} N$ on $M_f$,'' we
mean convergence in distribution of the restricted point processes.
This is the basic assumption made in this paper.

We generally assume the boundedness condition
\renewcommand{\theequation}{B}
\begin{equation}\label{eq:3.1}
E \Bigl[\Bigl(\sup_k Y_k\Bigr)^+\Bigr]<\infty.
\end{equation}
\renewcommand{\theequation}{\arabic{section}.\arabic{equation}}
\setcounter{equation}{0}

Let $\mathcal A_t=\sigma(N(\cdot\cap[0,t]\times\overline\R\cap
M_f))$, $t\in[0,1]$, denote the relevant filtration of the point
process $N$. A stopping time for $N$ or $N$-stopping time is a mapping
$T\dvtx \Omega\to[0,1]$ with $\{T\le t\}\in\mathcal A_t$ for each $t\in
[0,1]$. Denote by
\[
\overline Y_T:=\sup\{Y_k \dvtx  1\le k\le N(M_f), T=\tau_k\}, \qquad    \sup
\varnothing:= -\infty,
\]
the reward w.r.t. stopping time $T$.

Let $v\dvtx \overline M_f\to\overline\R$ be a continuous transformation
of the points of $N$ such that
%
\begin{equation}\label{eq:3.2}
\left.
\begin{array}{l}
v(t,x) \le ax^+ +b   \ \forall (t,x)\in M_f, \mbox{ with real
constants } a,b\ge0,  \\
v(t,\cdot) \mbox{ is for each $t$ a monotonically nondecreasing
function,}\\
v(\cdot, x) \mbox{ is for each $x$ a monotonically nonincreasing function.}
\end{array}
\right\}
\end{equation}

%

Define $c:=f(1)$ and for any guarantee value $x\in[c,\infty)$ and
$t\in[0,1)$ the \textit{optimal stopping curve} $\hat u$ of the
transformed Poisson process by
%
\begin{eqnarray}\label{eq:3.3}
\hat u(t,x) & :=& \sup\{E[ v(T,\overline Y_T\vee x)] \dvtx  T>t \mbox{ is
an } N\mbox{-stopping time} \},
\nonumber
\\[-8pt]
\\[-8pt]
\hat u(1,x) & :=& v(1,x).
\nonumber
\end{eqnarray}

It will be shown in the following proposition that the treshold
stopping time corresponding to $\hat u$ is an optimal
stopping time for the Poisson process.
For the basic notions of stopping of point processes; see \autorKRx
{-00a}, respectively, \autorFR.
The following proposition is the analogue of Proposition \ref
{prop:2.1} for continuous time Poisson processes. It is essential for
the solution of the $m$-stopping problem of $N$.

\begin{propos}[(Optimal stopping times larger than  $S$)]%
\label{prop:3.1}
Let $N$ satisfy the boundedness condition (\textup{B}), let $v$ satisfy
condition  \refeq{eq:3.2} and assume the following \textit{separation
condition} for the optimal stopping boundary $\hat u$:%
%
\renewcommand{\theequation}{$\hat{\mathrm{S}}$}
\begin{equation}\label{eq:3.4}%
\hat u(t,c)>\hat f(t):=v(t,f(t)) \qquad\forall  t\in[0,1).
\end{equation}
\renewcommand{\theequation}{\arabic{section}.\arabic{equation}}
\setcounter{equation}{2}
Then:
\begin{longlist}[(b)]
\item[(a)]$\hat u$ is continuous on $[0,1]\times[c,\infty]$ and for all
$(t,x)\in[0,1]\times[c,\infty]$ holds
%
\begin{eqnarray}\label{eq:3.5}
\hat{u}(t,x) &=& E \bigl[ v \bigl(T(t,x),\overline{Y}_{T(t,x)}\vee x\bigr) \bigr] \nonumber
\\[-8pt]
\\[-8pt]
&=& E \bigl[ v\bigl(T(t,x),\overline{Y}_{T(t,x)}\vee c\bigr)\vee v(1,x) \bigr]
\nonumber
\end{eqnarray}
with the optimal stopping time
\[
T(t,x) := \inf\{\tau_k >t \dvtx  v(\tau_k,Y_k) > \hat{u} (\tau_k,x)\},
  \qquad  \inf\varnothing:= 1.
\]
$\hat{u} (\cdot,x)$ is for $x\in[c,\infty]$ the optimal stopping
curve of the transformed Poisson process $\hat{N} := \sum_{k} \delta
_{(\tau_k, v(\tau_k,Y_k))} \mbox{ in } M_{\hat{f}}$ for the
guarantee value $v(1,x)$.

\item[(b)]
Let $S$ be an $N$-stopping time, let $Z\ge c$ be real $\mathcal
{A}_S$-measurable with $EZ^+<\infty$ and $\mathcal T(S)$ the set of
all $N$-stopping times $T$ with $T>S$ on $\{S<1\}$ and $T=1$ on $\{S=1\}
$. Then $T(S,Z)\in\mathcal T(S)$ is optimal in the sense that
%
\begin{eqnarray}\label{eq:3.6}
E \bigl[v \bigl(T(S,Z),\overline{Y}_{T(S,Z)}\vee Z\bigr) \mid\mathcal{A}_S\bigr] & =&
\hat{u}(S,Z) \nonumber
\\[-8pt]
\\[-8pt]
&\ge& E [v(T,\overline{Y}_{T}\vee Z) \mid\mathcal{A}_S] \qquad P\mbox{-a.s.}
\nonumber
\end{eqnarray}
for all $T\in\mathcal T(S)$.
\end{longlist}
\end{propos}

\begin{pf}
(a)
The statement in (a) is proved by discretization.
Since $\hat f$ is continuous and $\hat u(\cdot, c)$ is right
continuous there exists a monotonically nonincreasing, continuous
function $\hat f_2\dvtx [0,1]\to[\hat c,\infty)$, $\hat c:=\hat
f(1)=v(1,c)$ such that $\hat f<\hat f_2<\hat u(\cdot,c)$ on $[0,1)$.
Thus, for $t<1$, the sets $[0,t]\times\overline\R\cap M_{\hat f_2}$
are compact in $M_{\hat f}$.

For $x\in[c,\infty)$, $n\in\N$ and $1\le i\le2$ define
\[
M^{n}_{ {i}/{2^n}}(x) := \sup_{\tau_k\in ( ({i-1})/{2},
 {i}/{2} ]} v(\tau_k,Y_k\vee x).
\]
Consider the filtration $\mathcal A^n = (\mathcal A_{
{i}/{2^n}})_{1\le i\le2^n}$. Then $M_{ {i}/{2^n}}^n (x)$ is
$\mathcal A_{ {i}/{2^n}}$ measurable and $\mathcal A_{
{i}/{2^n}}$, $\sigma(M_{ ({i+1})/{2^n}}^n(x))$ are independent. We
define $w_n\dvtx [0,1]\times[c,\infty)\to\overline\R$ by
%
\begin{eqnarray}\label{eq:3.7}
   w_n(t,x) &:=& \sup\{E[ M^{n}_{T}(x)] \dvtx  T>t \mbox{ an } \mathcal{A}^n
\mbox{-stopping time}\} \qquad\mbox{for } t\in[0,1),\hspace*{-32pt}\nonumber
\\[-8pt]
\\[-8pt]
   w_n(1,x) &:=& v(1,x).\hspace*{-32pt}
\nonumber
\end{eqnarray}

Then for $t\in[0,1)$ by Proposition \ref{prop:2.1}, we have
\[
w_n(t,x) = E\bigl[ M^{n}_{T_n(t,x)}(x)\bigr] = V^n_{\lfloor2^nt\rfloor}(x)
\]
with the optimal $\mathcal{A}^{n}$-stopping time
\[
T_n(t,x) := \min \biggl\{t< \frac{i}{2^n}\le1 \dvtx
M^{n}_{ {i}/{2^n}}(x) > w_n\biggl(\frac{i}{2^n},x\biggr)  \biggr\},
\qquad \min\varnothing:= 1,
\]
and
%
\begin{eqnarray}\label{eq:3.8}
V^n_{2^n}(x) &:=& v(1,x),\nonumber
\\[-8pt]
\\[-8pt]
V^n_i(x) &:=& E \bigl[M^n_{ ({i+1})/{2^n}}(x)\vee V^n_{i+1}(x)\bigr], \qquad
i=2^n-1,\ldots , 0.
\nonumber
\end{eqnarray}
The function $w_n(\cdot,x)$ is monotonically nonincreasing and
constant on the intervals $[0, \frac{1}{2^n}), [\frac{1}{2^n},\frac
{2}{2^n}),\ldots ,[\frac{2^n-1}{2^n},1)$. We also have
\begin{eqnarray*}
 (1) &&\hspace*{5pt}    w_n(t,x)\ge\hat{u}(t,x)\qquad  \forall
t\in[0,1],\\
 (2) &&\hspace*{5pt}    w_n(t,x) \ge w_{n+1}(t,x)  \qquad  \forall  t\in[0,1].
\end{eqnarray*}

For the proof of  (1)  note that for any stopping time $T>t$,
$T_n:=\frac{\lceil{T2^n}\rceil}{2^n}$ is an $\mathcal
{A}^n$-stopping time with $T_n>t$ and $T_n-\frac{1}{2^n}<T\le T_n$. Therefore,
%
\begin{equation}\label{eq:3.9}
M^{n}_{T_n}(x) =\sup_{\tau_k\in   (T_n- {1}/{2^n},
T_n ]}v(\tau_k,Y_k\vee x) \ge v(T,\overline Y_T\vee x).
\end{equation}
This implies $w_n(t,x)\ge\sup\{E[ v(T,\overline Y_T\vee x)] \dvtx  T>t
\ N \mbox{-stopping time}\} = \hat{u}(t,x)$.

The proof of  (2)  is similar. If $T>t$ is an $\mathcal
{A}^{n+1}$-stopping time,\vspace*{-1pt} then $T':= \frac{\lceil{T2^n}\rceil}{2^n}$
is an $\mathcal{A}^n$-stopping time with $T'>t$ and $T'-\frac
{1}{2^n}< T \le T'$. Thus, as above, we obtain $w_n(t,x) \ge w_{n+1}(t,x)$.

Relations (1)  and  (2)  imply the existence of a monotonically
nonincreasing function $w(\cdot, x)\dvtx [0,1] \to\R\cup\{-\infty\}$
with $w(\cdot,x)\ge\hat u(\cdot,x)$ and $w_n(\cdot,x)\downarrow
w(\cdot,x)$ pointwise. It can be shown by our assumptions on $v$ and
$N$ that $w$ is continuous [see \autorF].

For $\omega\in\Omega$ with $\hat N(\omega,K)<\infty$ for all
compact $K\subset M_f$ and for $(t,x)\in[0,1]\times[c,\infty]$ and
$t_n\downarrow t$, we have the convergence
%
\begin{equation}\label{eq:3.10}
M^{n}_{T_n(t_n,x)}(x) \to v\bigl(T(t,x),\overline{Y}_{T(t,x)}\vee x\bigr)
\end{equation}
with the stopping time
%
\begin{eqnarray}\label{eq:3.11}
T(t,x) & :=& \inf\{\tau_k >t \dvtx  v(\tau_k,Y_k\vee x) > w(\tau_k,x)\}
\nonumber
\\[-8pt]
\\[-8pt]
& \stackrel{(*)}{=}& \inf\{\tau_k >t \dvtx  v(\tau_k,Y_k) >
w(\tau_k,x)\}, \qquad\inf\varnothing:= 1.
\nonumber
\end{eqnarray}

For the proof, note that monotone convergence of $w_n(\cdot, x)$ and
continuity of the limit $\omega$ implies uniform convergence from
above. Thus, for $x\in[c,\infty)$ points of $N$ on the graph of
$w(\cdot, x)$ are ignored by all stopping times $T_n(t,x)$ and
$T(t,x)$. The second equality $(*)$ holds since $w(t,x)\ge\hat
u(t,x)\ge v(t,x)$ and since by assumption $v(t,\cdot)$ is strictly
monotonically increasing. This implies by Fatou's lemma the following
sequence of inequalities:
\begin{eqnarray*}
\hat{u}(t,x) &\le& w(t,x) %
= \lim_{n\to\infty} w_n(t,x) = \lim_{n\to\infty}
E\bigl[M^{n}_{T_n(t,x)}(x)\bigr]\\
&\le& E \bigl[v\bigl(T(t,x),\overline{Y}_{T(t,x)}\vee x\bigr)\bigr] \le\hat{u}(t,x).
\end{eqnarray*}

Thus, $\hat u(\cdot, x)=w(\cdot,x)$ is continuous and $\hat u(t,x)=E[
v(T(t,x), \overline Y_{T(t,x)}\vee x)]$. As $w(t,x)\ge v(t,x)$ implies
that $\overline Y_{T(t,x)}> x$ for $T(t,x)<1$, we have $\hat
u(t,x)=E[v(T(t,x),\overline Y_{T(t,x)}\vee c)\vee v(1,x)]$, which means
that $\hat u(\cdot,x)$ is the optimal stopping curve of the Poisson
process $\hat N$ with guarantee value $v(1,x)$.

(b)
To prove optimality of the stopping time $T(S,Z)$, set $S_n:=\frac
{\lceil S2^n\rceil}{2^n}$. Then $S_n$ is an $\mathcal A^n$-stopping
time and by \eqref{eq:3.10} holds
%
\begin{equation}\label{eq:3.12}
M^{n}_{T_n(S_n,Z)}(Z) \to v\bigl(T(S,Z),\overline{Y}_{T(S,Z)}\vee Z\bigr)\qquad
P\mbox{-a.s.}
\end{equation}

Let $\mathcal T(S_n)$ be the set of all $\mathcal A^n$-stopping times
$T_n$ with $T_n>S_n$ on $\{S_n<1\}$ and $T_n=S_n$ on $\{S_n=1\}$. Let
$T\in\mathcal T(S)$. By discretization $T>S$ in general does not imply
$\frac{\lceil T2^n\rceil}{2^n} > \frac{\lceil S2^n\rceil}{2^n}$.
Thus, we modify the discretization and define $T_n := \frac{\lceil
{T2^n}\rceil}{2^n} \chi_{\{ {\lceil{T2^n}\rceil}/{2^n} > S_n\}}
+1\chi_{\{ {\lceil{T2^n}\rceil}/{2^n} = S_n\}} \in\mathcal
T(S_n)$. Then analogously to \refeq{eq:3.9}
\[
v(T,\overline{Y}_{T}\vee Z) \le M^{n}_{T_n}(Z) \chi_{\{ {\lceil
{T2^n}\rceil}/{2^n} > S_n\}} + v(T,\overline{Y}_{T}\vee Z) \chi_{\{
 {\lceil{T2^n}\rceil}/{2^n} = S_n\}}.
\]
This implies the inequalitites
\begin{eqnarray*}
&&E  [ v(T,\overline{Y}_{T}\vee Z) \mid\mathcal{A}_{S_n}
 ]
\\
&& \qquad \le E  [ M^{n}_{T_n}(Z) \mid\mathcal{A}_{S_n}  ] \chi_{\{
 {\lceil{T2^n}\rceil}/{2^n} > S_n\}}
\\
&&  \qquad  \quad   {} + E  [ v(T,\overline{Y}_{T}\vee Z)
\mid\mathcal{A}_{S_n}  ] \chi_{\{ {\lceil{T2^n}\rceil
}/{2^n} = S_n\}}
\\
&& \qquad \stackrel{(*)}{\le} \underbrace{E  \bigl[M^{n}_{T_n(S_n,Z)}(Z) \mid
\mathcal{A}_{S_n} \bigr]}_{= w_n(S_n,Z)} \chi_{\{ {\lceil
{T2^n}\rceil}/{2^n} > S_n\}}
\\
&&  \qquad  \quad   {} + E [v(T,\overline{Y}_{T}\vee Z)
\mid\mathcal{A}_{S_n} ] \chi_{\{ {\lceil{T2^n}\rceil
}/{2^n} = S_n\}}.
\end{eqnarray*}
$(*)$ holds by Remark \ref{rem:2.1}. Since we have\vspace*{-1pt} $M^n_{
{i}/{2^n}}(Z) = h(Y_i,Z)$, where $Y_i := N(\cdot\cap({\frac
{i-1}{2^n}},{\frac{i}{2^n}}]\times\overline\R\cap M_f)$ and with
$h\dvtx \mathcal{N}_R(M_f)\times[c,\infty)\to\overline\R$, $h(\sum_k
\delta_{(t_k,y_k)},x):=\sup_k v(t_k,y_k\vee x)$.

As $\mathcal{A}_S\subset\mathcal{A}_{S_n}$ we conclude
\begin{eqnarray*}
 &&E [ v(T,\overline{Y}_{T}\vee Z) \mid\mathcal
{A}_{S} ]
 \\
&& \qquad \le E  \bigl[ M^{n}_{T_n(S_n,Z)}(Z) \chi_{\{ {\lceil
{T2^n}\rceil}/{2^n} > S_n\}} \mid\mathcal{A}_{S} \bigr]
\\
&&  \qquad  \quad   {}
+ E  \bigl[ v(T,\overline{Y}_{T}\vee Z) \chi_{\{ {\lceil
{T2^n}\rceil}/{2^n} = S_n\}}  \mid  \mathcal
{A}_{S} \bigr]
\\
&&  \qquad = w_n(S_n,Z) E \bigl[\chi_{\{ {\lceil{T2^n}\rceil}/{2^n} > S_n\}
}  \mid  \mathcal{A}_{S} \bigr]
+ E  \bigl[ v(T,\overline{Y}_{T}\vee Z) \chi_{\{{\lceil
{T2^n}\rceil}/{2^n} = S_n\}}  \mid  \mathcal
{A}_{S} \bigr],
\end{eqnarray*}
and by the Lemma of Fatou we have by \refeq{eq:3.12}
\[
E  [ v(T,\overline{Y}_{T}\vee Z) \mid\mathcal{A}_{S} ] %
\le E \bigl[ v\bigl(T(S,Z),\overline{Y}_{T(S,Z)}\vee Z\bigr) \mid\mathcal
{A}_{S} \bigr]
= \hat{u}(S,Z).
\]
As $T>S$ was chosen arbitrary this implies (b).
\end{pf}

In the sequel, we need the following \textit{differentiability
condition} to be fulfilled.

(D)
 Assume that there is a version of the density
$g$ of $\mu$ on $M_f$ such that the intensity function
\[
G(t,y)=\int_y^\infty g(t,z)\,dz
\]
is continuous on $M_f\cap[0,1]\times\R$. Furthermore, we assume that
$\lim_{y\to\infty} y G(t,\break y)=0$ for all $t\in[0,1]$.

The following proposition determines the intensity function of
transformed Poisson processes.
\begin{propos}[(Intensity function of transformed
Poisson processes)] 
\label{prop:3.2}
Let $N=\sum\delta_{(\tau_k,Y_k)}$ be a Poisson process with
intensity function $G$ satisfying the boundedness condition  (\textup{B}). Let
$v\dvtx \overline M_f\to\overline\R$, $v=v(t,x)$ be a $C^1$-function
monotonically nonincreasing in $t$ and monotonically nondecreasing in
$x$ with $v(t,\infty)=\infty$ for all $t\in[0,1]$. Define $R(t,x) :=
(t,v(t,x))$ and $f_v(t):=v(t,f(t))$. Then $R(M_f)=M_{f_v}$,
$R^{-1}(t,y)=(t,\xi(t,y))$ with a $C^1$-function $\xi\dvtx M_{f_v}\to
\overline\R$.

$\widehat N:=\sum_k \delta_{(\tau_k,v(\tau_k,Y_k))}$ is a Poisson
process on $M_{f_v}$ with intensity measure $\widehat\mu=\mu\circ
R^{-1}$ and intensity fuction $\widehat G(t,y):=G(t,\xi(t,y))$,
$(t,y)\in M_{f_v}$.
\end{propos}

\begin{pf}
By \citeauthor{Resnick-87} [(\citeyear{Resnick-87}), Proposition 3.7], $\widehat N$ is a Poisson process
with intensity measure $\widehat\mu=\mu\circ R^{-1}$. The
transformation formula implies that the density $\widehat g$ of
$\widehat\mu$ is given by
\begin{eqnarray*}
\hat{g}(t,y) &=& g(R^{-1}(t,y)) |{\det J(R^{-1})(t,y)} |
 \\
&=& g(t,\xi(t,y)) \,\frac{\partial}{\partial y} \xi(t,y) = -\frac
{\partial}{\partial y} G(t,\xi(t,y)).
\end{eqnarray*}
\upqed
\end{pf}

After this preparation, we now consider the $m$-stopping problem for
Poisson processes. The aim is to solve
%
\begin{equation}\label{eq:3.13}
E [\overline{Y}_{T_1}\vee\cdots \vee\overline{Y}_{T_m}] = \mathrm{sup},
\end{equation}
where the supremum is over all $N$-stopping times\footnote{ $T_1<\cdots
<T_m\le1$ signifies that $T_{i-1}<T_i$ for each $i$ on $\{ T_{i-1}<1\}
$ and $T_i=1$ on $\{ T_{i-1}=1\}$.} $0\le T_1 < \cdots < T_m \le1$.

This problem has been considered for Poisson processes on $[0,1]\times
(c,\infty)$ already in \citet{SaarioSakaguchi-92} in the special case
of intensity functions of the form
%
\begin{equation}\label{eq:3.14}
G(t,y)=\lambda\bigl(1-F(y)\bigr)
\end{equation}
with $\lambda>0$ and $F$ a continuous distribution function with
$F(c)=0$. Equation \eqref{eq:3.14} models the case of i.i.d. random variables
arriving at Poisson distributed arrival times. \citet
{SaarioSakaguchi-95} derive for this case differential equations for
the optimal stopping curves. Explicit solutions are however not given
in any case. In the following, we extend these results to the case of
general intensities. We subsequently also identify classes of examples
of intensity functions which allow essentially explicit solutions.

In order to guarantee the existence of optimal $m$-stopping times, we
restrict ourselves in the following to the case where the lower
boundary is constant, $f\equiv c$. Define optimal $m$-stopping curves
for guarantee value $x\in[c,\infty)$, $m\in\N$, and $t\in[0,1)$ by$^2$
%
\begin{eqnarray}\label{eq:3.15}
u^m(t,x)&:=& \sup \{E [\overline{Y}_{T_1}\vee\cdots \vee\overline
{Y}_{T_m}\vee x] \dvtx   t<T_1<\cdots <T_m\le1\nonumber\\
&&\hspace*{150pt}  N\mbox{-stopping
times}  \},
\\
u^m(1,x)&:=& x.\nonumber
\end{eqnarray}
Further let $u^0(t,x):=x$ for $(t,x)\in[0,1]\times[c,\infty]$ and
$u^m(t) := u^m(t,c)$ for $t\in[0,1]$.

$u^m(\cdot,x)$ is called \textit{optimal $m$-stopping curve} of $N$ for
guarantee value $x$. Define the inverse function $\xi^m\dvtx \overline
M_{u^m}\to\overline R$ by
%
\begin{equation}\label{eq:3.16}
\xi^m(t,u^m(t,x))=x \qquad\mbox{for }  (t,x)\in[0,1]\times
[c,\infty].
\end{equation}
Further define $\gamma^m\dvtx [0,1]\times[c,\infty]\to\overline\R$ by
%
\begin{equation}\label{eq:3.17}
\gamma^m(t,x):=\xi^{m-1} (t,u^m(t,x))
\end{equation}
as well as
%
\begin{equation}\label{eq:3.18}
\gamma^m (t):=\gamma^m(t,c)=\xi^{m-1}(t,u^m(t)).
\end{equation}
Then $\gamma^m(t,x)>x$ iff $u^m(t,x)>u^{m-1}(t,x)$ and further
\[
y>\gamma^m(t,x)  \quad \Leftrightarrow \quad  u^{m-1}(t,y)>u^m(t,x).
\]

The optimal $m$-stopping for Poisson processes can be reduced by the
previous structural results to $m$ 1-stopping problem for transformed
Poisson processes. The transformations are given by the optimal
stopping curves $u^m$ or equivalently by the inverses $\gamma^m$---both sequences of curves are defined recursively. Thus, we consider the
transformed Poisson processes
%
\begin{equation}\label{eq:3.19}
N^m := \sum_k  \delta_{(\tau_k,u^{m-1}(\tau_k,Y_k))} \qquad\mbox{on } M_{u^{m-1}}.
\end{equation}
Define the (optimal) stopping times $T_\ell^m(t,k)$ with guarantee
value $x$ by
%
\begin{eqnarray}\label{eq:3.20}
 \quad T_1^m(t,x) &:=& \inf \{\tau_k>t \dvtx  Y_k> \gamma^{m}(\tau_k,x)
\},\nonumber
\\[-8pt]
\\[-8pt]
 \quad T_\ell^m(t,x) &:=& \inf  \bigl\{\tau_k>T_{\ell-1}^m(t,x)\dvtx  Y_k> \gamma
^{m-\ell+1}\bigl(\tau_k,\overline{Y}_{T_{\ell-1}^m(t,x)} \vee x\bigr) \bigr\}.
\nonumber
\end{eqnarray}
%
The following theorem characterizes the optimal stopping time as
threshold stopping time based on the optimal stopping curves.
These are given by a system of $m$ differential equations of first order.

\begin{theorem}[(Optimal  $  m$-stopping of Poisson processes)]%
\label{theo:3.3}
Let $f\equiv c$ and $N$ satisfy the boundedness condition (\textup{B}) and the
separation condition (\textup{S}), that is, $u^1(t)>c$ for $t\in[0,1)$. Let
$t_0(x):=\inf\{t\in[0,1] \dvtx  \mu((t,1]\times(x,\infty])=0 \}$.
\begin{longlist}[(b)]
\item[(a)]
Then for $m\in\N$, $(t,x)\in[0,1)\times[c,\infty)$ holds
\begin{eqnarray*}
u^m(t,x) &=& E  \bigl[ \overline{Y}_{T^{m}_1(t,x)}\vee\cdots \vee
\overline{Y}_{T^{m}_m(t,x)}\vee x  \bigr]
\\
&=& E  \bigl[ u^{m-1}\bigl(T^{m}_1(t,x), \overline{Y}_{T^{m}_1(t,x)}\vee x\bigr)
 \bigr]
\end{eqnarray*}
with optimal stopping times $(T_\ell^m(t,x))_{1\le\ell\le m}$
defined in \eqref{eq:3.20}.
\item[(b)]%
For $(t,x)\in A:=\{(t,x)\in(0,1]\times[c,\infty)\dvtx t<t_0(x)\}$ holds
$u^m(t,x)>u^{m-1}(t,x)$ while $u^m(t,x)=u^{m-1}(t,x)=x$ else. In
particular, $u^m(t)>u^{m-1}(t)$ for $t\in[0,1)$ and $u^m(\cdot,x)$ is
the optimal stopping curve of the transformed Poisson process $N^m$.
\item[(c)]%
Under the differentiability condition, (\textup{D}) $u^m(\cdot,x)$ solves the
differential equation
%
\begin{eqnarray}\label{eq:3.21}
\frac{\partial}{\partial t}  u^m(t,x) %
& =& -\int_{u^m(t,x)}^\infty G(t,\xi^{m-1}(t,y))\,dy, \qquad t\in
[0,1),\nonumber
\\[-8pt]
\\[-8pt]
u^m(1,x) &=& x.
\nonumber
\end{eqnarray}
\item[(d)]
For $x>-\infty$, \eqref{eq:3.21} has a unique solution. If $c=-\infty
$ and if
%
\begin{equation}\label{eq:3.22}
\liminf_{s\uparrow1} \frac{u(s)}{b(s)} < \infty,
\end{equation}
where $b(s):=E[\sup_{\tau_k>s}Y_k]$, then also in this case
$u^m=u^m(\cdot,-\infty)$ for $m\ge2$ is uniquely determined by
\eqref{eq:3.21}.
\end{longlist}
\end{theorem}

\begin{pf}
The proof is by induction in $m$. Our induction hypothesis is that the
statement of Theorem \ref{theo:3.3} holds and moreover that for any
$n$-stopping time $S$ and any $\mathcal A_S$-measurable $Z\ge c$ with
$EZ^+<\infty$ we have $P$-a.s.
\[
E \bigl[ Z\vee\overline{Y}_{T_1^m(S,Z)}\vee\cdots \vee\overline
{Y}_{T_m^m(S,Z)}  \mid \mathcal{A}_S\bigr] = u^m(S,Z) %
\ge E [ Z\vee\overline{Y}_{T_1}\vee\cdots \vee\overline{Y}_{T_m}
 \mid \mathcal{A}_S]
\]
for all $N$-stopping times $S<T_1<\cdots <T_m\le1$. Further,
%
\begin{equation}\label{eq:3.23a}
A = \{(t,x)\in[0,1]\times[c,\infty)\dvtx  u^m(t,x)>u^{m-1}(t,x)\}.
\end{equation}
For the one-stopping problem $m=1$ the statement of Theorem \ref
{theo:3.3} is contained in \autorFR. Proposition \ref{prop:3.1} with
$v(t,x):=x$ implies the first part of the induction hypothesis while
the second part follows from \autorFRx{Lemma 2.1(c)}.

For the induction step $m\to m+1$, we obtain for all stopping times
$S<T_1<T_2<\cdots <T_{m+1}\le1$ and $Z\ge c$ $\mathcal A_S$-measurable
by the induction hypothesis (note that $\mathcal A_S\subset\mathcal A_{T_1}$):
%
\begin{eqnarray}\label{eq:3.23b}
 &&E[ (Z\vee\overline{Y}_{T_1}) \vee\overline{Y}_{T_2} \vee
\cdots \vee\overline{Y}_{T_{m+1}} \mid \mathcal{A}_S]
\nonumber\\
&& \qquad \le E \bigl[ (Z\vee\overline{Y}_{T_1}) \vee\overline
{Y}_{T_1^{m}(T_1,Z\vee\overline{Y}_{T_1})} \vee\cdots \vee\overline
{Y}_{T_{m}^{m}(T_1,Z\vee\overline{Y}_{T_1})} \mid \mathcal{A}_S\bigr]
\\
&& \qquad = E[ u^m(T_1,Z\vee\overline{Y}_{T_1}) \mid\mathcal
{A}_S].\nonumber
\end{eqnarray}
This expression is maximized by Proposition \ref{prop:3.1} by
$T_1=T_1^{m+1}(S,Z)$ where
\[
T_1^{m+1}(t,x) := \inf\{\tau_k >t \dvtx  u^m(\tau_k,Y_k) > \hat{u}(\tau
_k,x)\}, \qquad \inf\varnothing:= 1.
\]
The maximizing value is given by $\hat{u}(S,Z)$.

For the proof, we need to show that $\hat{u} (t,c) > u^m(t)$ for $t\in
[0,1)$. We next establish this and at the same time show (\ref
{eq:3.23a}) for $m+1$.

Note that for $x\in[c,\infty)$
\begin{eqnarray*}
\hat{u}(t,x) & = & \sup \{E[u^m (T,\overline{Y}_T\vee x)] \dvtx  T>t
 \ N\mbox{-stopping time}  \}%
\\
& \ge& E \bigl[u^m\bigl(T_1^m(t,x),\overline{Y}_{T_1^m(t,x)}\vee
x\bigr)\bigr]%
\\
& \stackrel{(*)}{\ge}& E\bigl [u^{m-1}\bigl(T_1^m(t,x),\overline
{Y}_{T_1^m(t,x)}\vee x\bigr)\bigr]%
\\
& =& u^m(t,x)  \qquad\mbox{by induction hypothesis}.
\end{eqnarray*}
By \eqref{eq:3.23a}, we have strict inequality in $(*)$ if and only if
$P((T_1^m(t,x),\overline Y_{T_1^m(t,x)}) \in A)>0$. Using Lemma 2.4
in \autorFR, we see that this is equivalent to $\mu(A\cap M_{\gamma
^m(\cdot,x)}\cap(t,1]\times\R)>0$. This in turn is equivalent to
%
\begin{equation}\label{eq:3.24}
A\cap M_{\gamma^m(\cdot,x)}\cap(t,1]\times\R\not=\varnothing
\end{equation}
[since $\gamma^m(\cdot,x)$ is monotonically nonincreasing and by
definition of $A$]. We are going to show that this is fulfilled for all
points $(t,x)\in A$.

So let $(t,x)\in A$ and thus by induction hypothesis
$u^m(t,x)>u^{m-1}(t,x)$ or equivalently $\gamma^m(t,x)>x$. Under the
assumption that $M_{\gamma^m(\cdot,x)}\cap(t,1]\times\R\subset
A^c$, we obtain that also $(t,\gamma^m(t,x))\in A^c$ since $A^c$ is
closed. This implies that
\[
u^m(t,\gamma^m(t,x))=u^{m-1}(t,\gamma^m(t,x))=u^m(t,x).
\]
Since $u^m(t,\cdot)$ is strictly increasing, it follows that $\gamma
^m(t,x)=x$, which is a contradiction. Thus, \eqref{eq:3.24} holds true.

With the choice $S:=t$, $Z:=x$ further, we obtain
\begin{eqnarray*}
\hat u(t,x) &=& E\bigl[ u^m \bigl(T_1^{m+1} (t,x), \overline
Y_{T_1^{m+1}(t,x)}\vee x\bigr)\bigr] = u^{m+1}(t,x).
\end{eqnarray*}
Finally, in \eqref{eq:3.23b} holds
\[
T_l^m  \bigl(T_1^{m+1} (S,Z), Z\vee\overline Y_{T_1^{m+1}(S,Z)} \bigr) =
T_{l+1}^{m+1} (S,Z).
\]
By Proposition \ref{prop:3.1} $u^{m+1} (\cdot,x)$ is the optimal
stopping curve of the Poisson process $N^{m+1}=\sum_k \delta_{(\tau
_k,u^m(\tau_k,Y_k))}$ on $M_{u^m}$ at the guarantee value $x$. We
already proved that the separation condition is fulfilled for the
stopping of $N^{m+1}$ and by Proposition~\ref{prop:3.2} $N^{m+1}$ has
the intensity function $ G^{m+1}(t,y):=G(t,\xi^m(t,y))$.
The existence and uniqueness results for the differential equation
\eqref{eq:3.21} therefore follow with our assumption from the
corresponding result in \autorFR\ for the case $m=1$.
\end{pf}


\section{Explicit calculation of optimal $ m$-stopping curves}
\label{sec:4}

For the case of one-stopping problems, some classes of intensity
functions $G(t,y)$ have been introduced in \autorFR\ which allow to
determine optimal stopping curves in explicit form. Solving the
optimality equations in \eqref{eq:3.21} for the sequence of optimal
stopping curves for the $m$-stopping problem is in general much more
demanding. However, for some of the classes considered in \autorFR\
explicit solutions can be given also in the $m$-stopping case.

We consider intensity functions $G(t,y)$ of the form
%
\begin{equation}\label{eq:4.1}
G(t,y) = H  \biggl(\frac{y}{v(t)} \biggr)\frac{|v'(t)|}{v(t)}
\end{equation}
or
\begin{equation}
\label{eq:4.2}
  G(t,y)= H\bigl(y-v(t)\bigr)|v'(t)|
\end{equation}
as in \autorFR\ with $v(1)=0$ or $v(1)=\infty$ in case \eqref{eq:4.1}
and $v(1)=-\infty$ in case \eqref{eq:4.2}. For the general motivation
of these classes and these conditions, we refer to \autorFR. In
particular, we will see that the main application considered in this
paper to $m$-stopping of i.i.d. sequences
with discount and observation costs
is covered by these classes.

We first state the results in the three cases mentioned and then give
the proof.

\textit{Case 1:
 $G$ satisfies \eqref{eq:4.1} with $v$ monotonically
nonincreasing, $v(1)=0$.}
Here $c=0$. $H\dvtx (0,\infty]\to[0,\infty)$ is monotonically
nonincreasing continuous, $\int_0^\infty H(x)\,dx>0$ and we assume that
$v\dvtx [0,1]\to[0,\infty)$ is a $C^1$-function with $v>0$ on $[0,1)$.

We define
%
\begin{equation}\label{eq:4.3}
R^1(x):=x-\int_x^\infty H(y)\,dy,\qquad x\in(0,\infty),
\end{equation}
and assume that there exists some $r>0$ with $R^1(r)=0$.
%
%
%
%
Define $r_0:=0$, \mbox{$\Phi^0(x):=x$}.
Then for $m\ge1$ by induction holds:

The function $R^m\dvtx (r_{m-1},\infty)\to\R$ given by
%
\begin{equation}\label{eq:4.5}
R^m(x) := x-\int_x^\infty H (\Phi^{m-1}(y))\,dy
\end{equation}
has exactly one zero $r_m\in(r_{m-1},\infty)$ and the optimal
$m$-stopping curves are given for $(t,x) \in[0,1) \times[0,\infty]$ by
%
\begin{equation}\label{eq:4.6}
u^m(t,x) = \phi^m \biggl(\frac{x}{v(t)} \biggr)v(t),
\end{equation}
where $\phi^m \dvtx  [0,\infty]\to[r_m,\infty]$ is the inverse function
of $\Phi^m \dvtx  [r_m,\infty] \to[0,\infty]$,
\[
\Phi^m(x) := x\exp \biggl(-\int_x^\infty \biggl(\frac{1}{R^m(y)}
-\frac{1}{y} \biggr)\,dy \biggr).
\]
The system of functions $(R^m,\Phi^m)$, respectively, $(u^m,\phi^m)$ is by
\eqref{eq:4.6} recursively defined. In particular, it holds that
%
\begin{equation}\label{eq:4.7}
u^m(t)=r_m v(t)
\end{equation}
and thus determination of the optimal stopping curves is reduced to
finding a zero point of $\R^m$.

\textit{Case 2: $G$ satisfies \eqref{eq:4.1} with $v$ monotonically
nondecreasing, $v(1)=\infty$.}
Here $c=-\infty$. $H\dvtx (-\infty,\infty]\to[0,\infty)$ is
monotonically nonincreasing continuous, $\int_{-\infty}^0 H(x)\,
dx>0$, $\int_0^\infty H(x)\,dx=0$ and $\int_y^0\frac{H(x)}{-x}\,
dx<\infty$ for $y<0$. Further, we assume that $v\dvtx [0,1]\to[0,\infty]$
is a $C^1$-function with $v<\infty$ on $[0,1)$.

We define
\[
R^1(x):=x+\int_x^\infty H(y)\,dy,\qquad x\in(-\infty,\infty) ,
\]
and assume that there exists some $r<0$ with $R^1(r)=0$. Define
$r_0:=-\infty$, $\Phi^0(x):=x$. Then for $m\ge1$ by induction holds:

The function $R^m\dvtx (r_{m-1},0)\to\R$ defined by
\[
R^m(x):=x+\int_x^0 H(\Phi^{m-1}(y))\,dy
\]
has exactly one zero $r_m\in(r_{m-1},0)$ and the optimal $m$-stopping
curves are given for $(t,x)\in[0,1)\times\overline\R$ by
%
\begin{equation}\label{eq:4.8a}
u^m(t,x) =
\cases{\displaystyle
 x  ,&\quad if   $x\ge0 $,\cr\displaystyle
 \phi^m \biggl (\frac{x}{v(t)} \biggr)v(t)  ,&\quad if   $x< 0$,
}
\end{equation}
where $\phi^m\dvtx [-\infty,0]\to[r_m,0]$ is the inverse of $\Phi^m \dvtx
[r_m,0] \to[-\infty,0]$,
\[
\Phi^m(x) := x\exp\biggl (\int_x^0  \biggl(\frac{1}{y}-\frac
{1}{R^m(y)}  \biggr)\,dy \biggr).
\]
In particular, $u^m(t)=r_m v(t)$.

\textit{Case 3: $G$ satisfies \eqref{eq:4.2} with $v$
monotonically nonincreasing $v(1)=-\infty$.}
Then $c=-\infty$. $H\dvtx (-\infty,\infty]\to[0,\infty)$ is
monotonically nonincreasing continuous, $\int_{-\infty}^\infty H(x)\,
dx>0$ and $\int_z^\infty\int_y^\infty H(x)\,dx\,dy<\infty$ for
$z\in\R$. Further, we assume that $v\dvtx [0,1]\to[-\infty,\infty)$ is
a $C^1$-function with $v>-\infty$ on $[0,1)$.

We define
\[
R^1(x):=1-\int_x^\infty H(y)\,dy,\qquad x\in\R,
\]
and assume that there exists some $r\in\R$ such that $R^1(r)=0$.
Define $r_0:=-\infty$, $\Phi^0(x):=x$. Then for $m\ge1$ by induction holds:

The function $R^m\dvtx (r_{m-1}, \infty)\to\R$ defined by
\[
R^m(x) := 1-\int_x^\infty H(\Phi^{m-1}(y)) \,dy
\]
has exactly one zero $r_m\in(r_{m-1},\infty)$. The optimal
$m$-stopping curves are given for $(t,x)\in[0,1) \times\overline\R$ by
%
\begin{equation}\label{eq:4.9a}
u^m(t,x) = \phi^m\bigl(x-v(t)\bigr)+v(t),
\end{equation}
where $\phi^m\dvtx \overline\R\to[r_m,\infty]$ is the inverse of $\Phi
^m\dvtx [r_m,\infty]\to\overline\R$,
\[
\Phi^m(x) := x-\int_x^\infty \biggl(\frac{1}{R^m(y)}-1  \biggr)\,dy.
\]
We have $u^m(t)=r_m+v(t)$.

\begin{pf}
We only give the proof of Case 2. The proof of both other cases is
similar. The proof is by induction in $m$ where we additionally include
that $R^m\ge R^{m-1}$ and thus $\Phi^m\ge\Phi^{m-1}$.

In the case $m=1$, the statement has been shown in \autorFR\ [with
$r_0:=-\infty$, $\Phi^0(x):=x$, $R^0(x):=x$].

Induction step $m\to m+1\dvtx  u^{m+1}(\cdot,x)$ is the optimal stopping
curve of $N^{m+1}$ at the guarantee value $x$. $N^{m+1}$ has the
intensity function
\[
G^{m+1}(t,y) = H \biggl (\Phi^m\biggl (\frac{y}{v(t)} \biggr) \biggr)
\frac{v'(t)}{v(t)} \qquad\mbox{for } (t,y) \in M_{u^m}.
\]
Thus, $G^{m+1}$ again is of type \eqref{eq:4.1} and we have to check
the conditions of Case~2 in \autorFR, who deal with optimal
one-stopping w.r.t. this type of intensity functions. First, we note
that $R^{m+1}$ has a zero in $(r_m,0)$ since $\Phi^m(x)\ge\Phi
^{m-1}(x)$ and thus $R^{m+1}\ge R^m$. Further by substitution, we have
\[
\int_y^0 \frac{H(\Phi^m(x))}{-x} \, dx \stackrel{\mathrm{Subst.}}{=}
\int_{\Phi^m(y)}^0\frac{H(z)}{-z} \frac{-z}{\phi^m(z)}  (\phi
^m)'(z)\, dz < \infty,
\]
as $\lim_{z\to0}\frac{-z}{\phi^m(z)}=1$ and $\lim
_{z\to0}(\phi^m)'(z)=1$. Thus, the conditions\vspace*{1pt} hold true and the
result follows.
\end{pf}

For intensity functions $G$ not of the form as in \eqref{eq:4.1},
\eqref{eq:4.2} the optimality differential equations in Theorem \ref
{theo:3.3} typically can only be solved numerically. In some cases,
however, one can derive bounds for the optimal stopping curves $u^m(t,x)$
which can be used to derive necessary uniform integrability and separation
conditions [see \autorFx{pages 60--62}] for the following approximation result.

\section{Approximation of $ m$-stopping problems}
\label{sec:5}

In this section, an extension of the approximation results in
\citeauthor{KuehneRueschendorf-04} [(\citeyear{KuehneRueschendorf-04}, Theorem 2.1]
and \citeauthor{FallerRueschendorf-09} [(\citeyear{FallerRueschendorf-09}), Theorem 4.1], for optimal
one-stopping problems for dependent sequences is given to the class of
$m$-stopping problems. For the special case of i.i.d. sequences with
distribution function $F$ in the domain of the Gumbel extreme value
distribution $\Lambda$ a corresponding approximation result was given
in the case $m=2$ in \autorKRx{-02}. The following result concerns the
dependent case and needs a new technique of proof which is based on
discretization. The main result of this section states that under some
conditions convergence of the finite imbedded point processes $N_n$ to
a Poisson process $N$ implies approximation of the stopping behavior.

We use the same general assumptions as in Section~4 of \autorFR\ as well
as the notation in Section~2 for the Poisson process $N$. In
particular, $\gamma^1, \ldots , \gamma^m$ are the functions defined in
\eqref{eq:3.17}. Further, the lower boundary curve $f$ of $N$ is given
by $f\equiv c$, $N$ is a Poisson process on $[0,1]\times(\overline\R
\setminus\{c\})$ and $\mathcal F^n$ are the canonical filtrations
induced by the imbedded point process $N_n$ and\vspace*{-1pt} we assume the
convergence condition
$N_n \stackrel{d}{\rightarrow} N$ on $M_f$ as throughout this paper [see (\ref
{eq:1.2}) and the introduction of Section~\ref{sec:3}].

The first result is an extension of Proposition 2.4 in \autorKRx{-00a}\
on the convergence of threshold stopping times to the case
$m\ge1$. For the technically involved proof, we refer to \autorFx{Lemma 2.6}.

\begin{propos}[(Convergence of multiple threshold stopping times)]%
\label{prop:5.1}
Let $(t,x)\in[0,1]\times[c,\infty)$ be fixed and let $v^m_n\dvtx [0,1]\to
\overline\R$ be functions such that $v^m_n \to\gamma^m (\cdot,x)$
uniformly on any interval $[0,s]$ with $s<1$. Define the corresponding
threshold stopping times
\begin{eqnarray*}
\hat{T}^{n,m}_1(t,x)
& :=& \min \biggl\{ tn < i \le n-m+1 \dvtx  X^n_i > v^m_n \biggl(
\frac{i}{n} \biggr) \biggr\},
\\
\hat{T}^{n,m}_\ell(t,x) & :=& \min \biggl\{ \hat{T}^{n,m}_{\ell
-1}(t,x)< i \le n-m+\ell\dvtx\\
&&\hphantom{\min \biggl\{}
X^n_i > \gamma^{m-\ell+1} \biggl(\frac{i}{n}, X^n_{\hat
{T}^{n,m}_{\ell-1}(t,x)}\vee x \biggr) \biggr\}
\end{eqnarray*}
for $2\le\ell\le m$. If $N_n \stackrel{d}{\to} N$ on $M_c$, we
obtain convergence
%
\begin{equation}\label{eq:5.1}
 \qquad  \biggl(\frac{\hat{T}^{n,m}_\ell(t,x)}{n}, X^n_{\hat
{T}^{n,m}_\ell(t,x)} \vee x \biggr)_{1\le \ell\le m} \stackrel
{d}{\longrightarrow} \bigl(T^m_\ell(t,x), \overline{Y}_{T^m_\ell(t,x)}
\vee x\bigr)_{1\le\ell\le m}.
\end{equation}
\end{propos}

Let now $W_k^{n,m}(x)$ be the stopping thresholds for the $m$ stopping
of $X_1^n, \ldots , X_n^n$ and the filtration $\mathcal F^n$ (see
Section \ref{sec:2}). The optimal \textit{$m$-stopping curves} w.r.t.
$\mathcal F^n$ are defined as follows. For $t\in[0,\frac{n-m+1}{n})$
and $x\in\overline\R$ let
\[
u^m_n(t,x) := W^{n,m}_{\lfloor tn\rfloor}(x)
\]
and $u^m_n(t,x):= W^{n,m}_{n-m+1}(x)$ for $t\in[\frac{n-m+1}{n},1]$.

More explicitly, we have for $t\in[0, \frac{n-m+1}{n})$ (see Theorem
\ref{theo:2.2})
%
\begin{eqnarray}\label{eq:5.2}
u^m_n(t,x) &=&
\esssup \bigl\{ E\bigl[ X^n_{T_1}\vee\cdots \vee X^n_{T_m}\vee x \mid
\mathcal{F}^{n}_{\lfloor tn\rfloor}\bigr] \dvtx  tn <T_1<\cdots <T_m\le n \nonumber\\
 &&\hspace*{196.5pt}\mathcal{F}^n\mbox{-stopping times}
 \bigr\}
 \\
&=& E\bigl[ X^n_{T^{n,m}_1(t,x)} \vee\cdots \vee X^n_{T^{n,m}_m(t,x)} \vee x
\mid\mathcal{F}^n_{\lfloor tn\rfloor}\bigr]\qquad P\mbox{-a.s.}
\nonumber
\end{eqnarray}
The corresponding optimal $m$-stopping times are given by
%
\begin{eqnarray}\label{eq:5.3}
 \quad T^{n,m}_1(t,x) %
&:=& \min \biggl\{ tn< i\le n-m+1 \dvtx  u^{m-1}_n\biggl(\frac{i}{n},X^n_i\biggr) >
u^m_n \biggl(\frac{i}{n},x\biggr)  \biggr\},\hspace*{-22pt}
\nonumber\\
 \quad T^{n,m}_\ell(t,x) %
&:=& \min \biggl\{ T^{n,m}_{\ell-1}(t,x) < i \le n-m+\ell\dvtx  \\
 \quad &&\hphantom{\min \biggl\{}
u^{m-\ell}_n \biggl(\frac{i}{n},X^n_i\biggr)> u^{m-\ell+1}_n\biggl(\frac{i}{n},
M^{n,m}_{\ell-1, i}\vee x\biggr)  \biggr\}\nonumber
\end{eqnarray}
for $2\le\ell\le m$, where $M^{n,m}_{j, i}:=
X^n_{T^{n,m}_{j}(t,x)}\chi_{\{T^{n,m}_j(t,x) \le  i\}}$.

$u_n^m(\cdot,x)$ is right continuous and a piecewise constant curve in
the space of random variables. We have the iterative representation
(see Theorem \ref{theo:2.2})
\[
u^m_n(t,x) = E \biggl[ u^{m-1}_n  \biggl(\frac{T^{n,m}_1
(t,x)}{n}, X^n_{T^{n,m}_1 (t,x)} \vee x \biggr)  \Big|  \mathcal
{F}^{n}_{\lfloor tn\rfloor} \biggr] \qquad P\mbox{-a.s.}
\]
Further, $u_n^m$ are monotone in the sense that for $0\le s\le t\le1$
\[
u^m_n(s,x)  \ge  E  \bigl[ u^m_n(t,x) \mid \mathcal
{F}^{n}_{\lfloor sn\rfloor}  \bigr] \qquad P\mbox{-a.s.}
\]

In the opposite direction, we obtain for $0\le s\le t\le1$
%
\begin{equation}\label{eq:5.4}
u^m_n(s,x) \le E  \biggl[ \max_{s< {i}/{n} \le t} u_n^{m-1}
\biggl(\frac{i}{n},X^n_i\biggr) \vee u_n^m(t,x) \mid\mathcal{F}^{n}_{\lfloor
sn\rfloor} \biggr] \qquad P\mbox{-a.s.}
\end{equation}

This follows inductively from the recursive definition of the
thresholds $W_\ell^m(x)$. We also need the following further
conditions [for motivation, see \autorFR]:

\begin{longlist}[(L)]
\item[(A)] \textit{Asymptotic independence
condition.}
For $0\le s<t\le1$
\[
P \Bigl (\max_{s< {i}/{n}\le t} X^n_i \le x \mid\mathcal{F}^{
n}_{ \lfloor sn \rfloor} \Bigr) \stackrel{P}{\longrightarrow} P
\Bigl(\sup_{s<\tau_k\le t} Y_k\le x \Bigr) \qquad\forall
x\in(c,\infty).
\]
\item[(U)]
 \textit{Uniform integrability
condition.}
$M_n^+$, with $M_n := \max_{1\le i\le n} X^n_i$, is uniformly
integrable and $E [ \limsup_{n\to\infty} M_n^+] < \infty.$
\item[(L)]
\textit{Uniform integrability from
below.}
For some sequence $(v_n)_{n\in\N}$ of monotonically nonincreasing
functions $v_n\dvtx [0,1]\to\R\cup\{-\infty\}$ with $v_n\to u$
pointwise, for all $t\in[0,1)$ and the corresponding threshold
stopping times
\[
\hat{T}_n(t) := \min\biggl\{ tn<i\le n \dvtx  X^n_i > v_n \biggl(\frac
{i}{n} \biggr) \biggr\}
\]
holds
%
\begin{equation}\label{eq:5.5a}
\lim_{s \uparrow1} \limsup_{n\to\infty} E \bigl[ X^{n}_{ \hat{T}_n(t)}
\chi_{\{\hat{T}_n(t)> sn\}}\bigr] = 0.
\end{equation}
\end{longlist}

A modified version of (L) is the condition (L$^m$):

\begin{longlist}[(L$^m$)]
\item[(L$^m$)]
For $m\in\N$, there exists some sequence of monotonically
nonincreasing functions $v_n^m\dvtx [0,1]\to\overline\R$ such that
$v_n^m\to\gamma^m(\cdot,-\infty)$ pointwise and further the
corresponding threshold stopping times
\[
\hat T_1^{n,m}(t) := \min\biggl\{tn<i\le n-m+1 \dvtx  X_i^n>v_n^m\biggl(\frac{i}{n}\biggr)\biggr\}
\]
satisfy
\[
\lim_{s \uparrow1} \limsup_{n\to\infty} E\bigl[ X^n_{\hat
{T}^{n,m}_{1}(t)} \chi_{\{\hat{T}^{n,m}_{1}(t) > sn \}}\bigr] = 0.
\]
\end{longlist}

Condition (L$^m$) in combination with (U) implies uniform integrability
of $(X_{\hat T_1^{n,m}(t)}^n)_{n\in\N}$. Denote
\[
T^{n,m}_\ell:= T^{n,m}_\ell(0,c) \quad\mbox{and} \quad T^{m}_\ell
:= T^{m}_\ell(0,c).
\]

\begin{theorem}[(Approximation of $  m$-stopping problems)]%
\label{theo:5.2}
  Assume that $N_n\stackrel{d}{\to}N$ on $[0,1]\times
(\overline\R\setminus\{c\})$ and also assume conditions (\textup{A}) and (\textup{U}).
In case $c=-\infty$ also assume the modified uniform integrability
condition (\textup{L}$^{{m}}$).

\begin{longlist}[(b)]
\item[(a)]%
For all $(t,x)\in[0,1]\times[c,\infty)$ holds
\[
u_n^m(t,x) \stackrel{P}{\longrightarrow} u^m(t,x).
\]
If $c\in\R$ assume $X^n_n\stackrel{L^1}{\to} c$. Then we have in particluar
%
\begin{equation}\label{eq:5.5b}
E  [X^n_{T^{n,m}_1} \vee\cdots \vee X^n_{T^{n,m}_m}  ] \to u^m(0).
\end{equation}
\item[(b)]%
In case $(X_i^n)_{1\le i\le n}$ are independent random variables and if\vspace*{-2pt}
for $c\in\R$ we assume that $\mu(M_{\gamma^m})=\infty$ or
$X_{n-i}^n \stackrel{P}{\longrightarrow} c$ for $i=0,\ldots ,m-1$,
then we obtain
\[
 \biggl(\frac{T^{n,m}_\ell}{n}, X^n_{T^{n,m}_\ell} \biggr)_{1\le
\ell\le m} \stackrel{d}{\to}
(T^m_\ell, \overline{Y}_{T^m_\ell} \vee c)_{1\le\ell\le m}.
\]
\item[(c)]%
If $c\in\R$ and $X^n_n\stackrel{L^1}{\to} c$, then
\begin{eqnarray}
\hat{T}^{n,m}_1 %
&:=& \min \biggl\{ 1\le i \le n-m+1 \dvtx  X^n_i > \gamma^m\biggl(\frac{i}{n},c\biggr)
 \biggr\},\nonumber
 \\
\hat{T}^{n,m}_\ell%
&:=& \min \biggl\{\hat{T}^{n,m}_{\ell-1} < i \le n-m+\ell\dvtx
X^n_i > \gamma^{m-\ell+1}\biggl(\frac{i}{n}, X^n_{\hat{T}^{n,m}_{\ell
-1}} \vee c\biggr) \biggr\},\nonumber
\\
\eqntext{2\le\ell\le m,}
\end{eqnarray}
defines an asymptotically optimal sequence of $m$-stopping times, that
is, convergence as in \eqref{eq:5.5b} holds for these stopping times.
In case $c=-\infty$,\vspace*{1pt}
\begin{eqnarray}
\hat{T}^{n,m}_1 %
&:=& \min \biggl\{1\le i \le n-m+1 \dvtx  X^n_i > v^m_n\biggl(\frac{i}{n}\biggr) \biggr\},
\nonumber\\
\hat{T}^{n,m}_\ell%
&:=& \min \biggl\{\hat{T}^{n,m}_{\ell-1} < i \le n-m+\ell\dvtx
X^n_i > \gamma^{m-\ell+1}\biggl(\frac{i}{n}, X^n_{\hat{T}^{n,m}_{\ell
-1}}\biggr) \biggr\},
\nonumber\\[1pt]
\eqntext{2\le\ell\le m,}
\end{eqnarray}
are asymptotically optimal stopping times, where $v_n^m$ are the
threshold functions from condition $(\mathrm{L}^m)$.
\end{longlist}
\end{theorem}

\begin{pf}
Since we use point process convergence on $[0,1]\times(\overline\R
\setminus\{c\})$ and canonical filtrations, we can apply the Skorohod
theorem and hence we assume w.l.o.g. $P$-a.s. convergence of the point
processes.

 (a)
Consider at first the case $c\in\R$. Let $t\in[0,1)$ be a fixed
element. We introduce at first discrete majorizing stopping problems.
For $m\ge1$ and $k>m$, define the discrete time points\vspace*{1pt}
\[
a_i^k :=  \biggl(1- \frac{i}{k} \biggr)t + \frac{i}{k} 1, \qquad0\le
i\le k,
\]
and discrete time random variables\vspace*{1pt}
\[
X^{n,k}_i := \max_{ {j}/{n}\in(a^k_{i-1}, a^k_i]} X^n_j\vee c
\qquad\mbox{for }1\le i\le k,
\]
and consider the filtration $\mathcal{F}^{n,k}:=(\mathcal
{F}^{n,k}_i)_{0\le  i\le  k}$ with $\mathcal{F}^{n,k}_i :=\mathcal
{F}^{n}_{\lfloor a^k_i n\rfloor}$. The corresponding $m$-stopping
curves are given inductively for $m\ge1$ by backward induction for
$i=k,\ldots ,0$ by\vspace*{1pt}
\begin{eqnarray}
\sideset{^m}{^{n,k}_{k-m+1}}{{W}}(x) &:= &x,
\nonumber\\
\sideset{^m}{^{n,k}_i}{{W}}(x) %
&:=& E [ \sideset{^{m-1}}{^{n,k}_{i+1}}{{W}}(X^{n,k}_{i+1})
\vee\sideset{^m}{^{n,k}_{i+1}}{{W}}(x) \mid
\mathcal{F}^{n,k}_i]
\nonumber\\
\eqntext{\mbox{for } i= k-m,\ldots ,0.}
\end{eqnarray}
These stopping problems majorize the original $m$-stopping problem,\vspace*{1pt}
\begin{eqnarray*}
\sideset{^m}{^{n,k}_{0}}{{W}}(x)
&=& \esssup\{ E[ X^{n,k}_{T'_1}\vee\cdots \vee X^{n,k}_{T'_m}\vee x \dvtx
\mathcal{F}^{n,k}_{0}] \dvtx  \\
&&\hphantom{\esssup\{}  0<T'_1<\cdots <T'_m\le k \ \mathcal{F}^{n,k}\mbox{-stopping times} \}
\\
&\stackrel{(*)}{=}& \esssup\{ E[ X^{n,k}_{T'_1}\vee\cdots \vee
X^{n,k}_{T'_m}\vee x \dvtx  \mathcal{F}^{n,k}_{0}] \dvtx  \\
&&\hphantom{\esssup\{}  0<T'_1\le\cdots \le T'_m\le k \ \mathcal
{F}^{n,k}\mbox{-stopping times} \}
\\
&\ge& \esssup\bigl\{ E\bigl[ X^n_{T_1}\vee\cdots \vee X^n_{T_m}\vee x \mid
\mathcal{F}^{n}_{\lfloor t n\rfloor}\bigr] \dvtx  \\
&& \hphantom{\esssup\{} tn<T_1<\cdots <T_m\le n \ \mathcal{F}^{n}\mbox
{-stopping times} \bigr\}
\\
&=& u^m_n(t,x)\qquad P\mbox{-a.s.},
\end{eqnarray*}
since for all $\mathcal{F}^{n}$-stopping times $tn<T_1<\cdots <T_m\le
n$ it holds that $T'_i := \lceil\frac{1}{1-t}(\frac
{T_i}{n}-t)k\rceil>0$ are $\mathcal{F}^{n,k}$-stopping times with
$a^k_{T'_i-1}<\frac{T_i}{n}\le a^k_{T'_i}$, thus $X^{n,k}_{T'_i} \ge
X^n_{T_i}$.
For the proof of $(*)$ define for $\mathcal{F}^{n,k}$-stopping times
$0<T'_1\le\cdots \le T'_m\le k$ the $\mathcal{F}^{n,k}$-stopping times
$0<T_1^*< \cdots < T_m^*\le k$ by
\begin{eqnarray*}
T_1^* && := T'_1\wedge(k-m+1),\\
T_\ell^* && :=  \bigl((T'_\ell+1) \chi_{\{T^*_{\ell-1} = T'_\ell\}
} +T'_\ell\chi_{\{T^*_{\ell-1}< T'_\ell\}} \bigr) \wedge(k-m+\ell
), \qquad\ell=2,\ldots ,m.
\end{eqnarray*}
We will prove convergence as $n\to\infty$ to the stopping problem of
\[
Y^k_i := \sup_{\tau_l\in(a_{i-1}^k, a_i^k]} Y_l\vee c \qquad\mbox
{for } 1\le i\le k,
\]
with filtrations $\mathcal{A}^{k} := (\mathcal{A}^{k}_i)_{1\le i\le
k}$, $\mathcal{A}^{k}_i := \mathcal{A}_{a_i^k}$ and optimal thresholds
\begin{eqnarray*}
\sideset{^m}{^k_{k-m+1}}{{u}} (x) &:= & x,\\
\sideset{^m}{^k_i}{{u}} (x) & := & E[ \sideset
{^{m-1}}{^k_{i+1}}{{u}} (Y^k_{i+1})\vee\sideset
{^m}{^k_{i+1}}{{u}} (x)] \qquad\mbox{for } i= k-m,\ldots ,0.
\end{eqnarray*}
By definition for $i\le k-m$ holds
\begin{eqnarray*}
\sideset{^m}{^k_i}{{u}} (x) &=& V\bigl(\sideset
{^{m-1}}{^k_{i+1}}{{u}} (Y^k_{i+1})\vee x,\ldots , \sideset
{^{m-1}}{^k_{k-m+1}}{{u}} (Y^k_{k-m+1})\vee x\bigr)\\
&=& \sup\{ E[ \sideset{^{m-1}}{^k_T}{{u}} (Y^k_T)\vee x] \dvtx
i<T\le k-m+1 \ \mathcal{A}^k\mbox{-stopping times}\}
\\
&=& \sideset{^m}{^k}{{u}} (a^k_i,x),
\end{eqnarray*}
where $\sideset{^m}{^k}{{u}} (\cdot,x)$ are the optimal
stopping curves of the processes
\[
\sideset{^m}{^k}{{N}}%
:= \sum_{i=1}^{k-m+1} %
\delta_{(a^k_i, \sideset{^{m-1}}{^k_i}{u} (Y^k_i))} =
\sum_{i=1}^{k-m+1} %
\delta_{(a^k_i, \sideset{^{m-1}}{^k}{u} (a^k_i,Y^k_i))}
\]
at guarantee value $x$.

At first we establish that for any $i$ the random variable $Y_{i+1}^k$
is\vspace*{-1pt} independent of the $\sigma$-algebra $\mathcal F_i^k := \sigma
(\bigcup_{n\in\N} \mathcal F_i^{n,k})$.

For the proof, note that by condition (A)
\[
P(X^{n,k}_{i+1} \in\cdot\mid\mathcal{F}^{n,k}_{i}) \stackrel
{P}{\longrightarrow} P(Y^k_{i+1}\in\cdot).
\]
Thus, we obtain by the continuous mapping theorem that for any
continuous $f\dvtx \overline\R\to[0,1]$ we have
\[
P\bigl(f(X^{n,k}_{i+1}) \in\cdot \mid\mathcal{F}^{n,k}_{i}\bigr) \stackrel
{P}{\longrightarrow} P\bigl(f(Y^k_{i+1})\in\cdot \bigr).
\]
This implies using uniform integrability that
\[
E [f(X^{n,k}_{i+1}) \mid \mathcal{F}^{n,k}_{i}] \stackrel
{L^1}{\longrightarrow} E [f(Y^k_{i+1})].
\]
On the other hand, by point process convergence it holds that\vspace*{-3pt}
$X^{n,k}_{i+1} \to Y^k_{i+1} P$-a.s. and thus also $f(X^{n,k}_{i+1})
\stackrel{L^1}{\longrightarrow} f(Y^k_{i+1})$. This implies
$L^1$-convergence of conditional expectations:
\[
E [f(X^{n,k}_{i+1}) \mid \mathcal{F}^{n,k}_{i}] \stackrel
{L^1}{\longrightarrow} E [f(Y^k_{i+1}) \mid\mathcal{F}^{k}_i].
\]
In consequence, we obtain $E [f(Y^k_{i+1})] = E [f(Y^k_{i+1}) \mid
\mathcal{F}^{k}_i]$ $P$-a.s. for all continuous functions $f\dvtx \overline
\R\to[0,1]$, and thus independence of $\mathcal{F}^{k}_i$ and
$\sigma(Y^k_{i+1})$.

The next point to establish is proved by induction in $m$. The
induction hypothesis is:
\begin{longlist}[(1)]
\item[(1)]%
For all $k>m$, $x\in[c,\infty)$ and $i=k-m+1,\ldots ,0$
\[
\sideset{^m}{^{n,k}_i}{{W}}(x) \stackrel{P}{\longrightarrow}
\sideset{^m}{^k_i}{{u}}(x), \qquad n\to\infty.
\]
\item[(2)]
For all $s\in[t,1]$ and all $x\in[c,\infty)$, we further have
\[
\sideset{^m}{^k}{{u}}(s,x) \to u^m(s,x), \qquad k\to\infty.
\]
\end{longlist}

We do the induction step for $m-1\to m$: Assertion (1) we shall prove
by backward induction on $i$: For $i=k-m+1$ the assertion is trivial.
We now consider the induction step from $i+1$ to $i$: From the
induction hypothesis, we know that
\[
\sideset{^{m-1}}{^{n,k}_{i+1}}{{W}}(x) \stackrel
{P}{\longrightarrow} \sideset{^{m-1}}{^k_{i+1}}{{u}}(x),
\qquad n\to\infty,
\]
for all $x\in[c,\infty)$.
From this, the monotonicity of $\sideset{^{m-1}}{^{n,k}_{i+1}}{
{W}}(x)$ in $x$ and the continuity of $\sideset
{^{m-1}}{^k_{i+1}}{{u}}(x)$ in $x$ we can conclude that
\[
\sideset{^{m-1}}{^{n,k}_{i+1}}{{W}}(X^{n,k}_{i+1}) \stackrel
{P}{\longrightarrow} \sideset{^{m-1}}{^k_{i+1}}{
{u}}(Y^k_{i+1}),\qquad n\to\infty.
\]
For details, see \autorF.
By the induction hypothesis for $i$, we also know that
\[
\sideset{^m}{^{n,k}_{i+1}}{{W}}(x) \stackrel
{P}{\longrightarrow} \sideset{^m}{^k_{i+1}}{{u}}(x),\qquad
n\to\infty,
\]
for $x\in[c,\infty)$, implying
\[
\sideset{^{m-1}}{^{n,k}_{i+1}}{{W}}(X^{n,k}_{i+1})\vee
\sideset{^m}{^{n,k}_{i+1}}{{W}}(x) \stackrel
{L^1}{\longrightarrow} \sideset{^{m-1}}{^k_{i+1}}{
{u}}(Y^k_{i+1})\vee\sideset{^m}{^k_{i+1}}{{u}}(x),\qquad n\to
\infty.
\]
From this, we get
\begin{eqnarray*}
 &&E[ \sideset{^{m-1}}{^{n,k}_{i+1}}{
{W}}(X^{n,k}_{i+1})\vee\sideset{^m}{^{n,k}_{i+1}}{{W}}(x)\mid
\mathcal{F}^{n,k}_i]   \\
&& \qquad \stackrel{L^1}{\longrightarrow} E[ \sideset
{^{m-1}}{^k_{i+1}}{{u}}(Y^k_{i+1})\vee\sideset
{^m}{^k_{i+1}}{{u}}(x)\mid\mathcal{F}^{k}_i]
\end{eqnarray*}
as $n\to\infty$. The expression on the left-hand side equals
$\sideset{^m}{^{n,k}_i}{{W}}(x)$, and since $\sigma
(Y^k_{i+1})$ and $\mathcal{F}^{k}_i$ are independent as shown above,
the right-hand side equals $\sideset{^m}{^k_i}{{u}}(x)$. This
completes the induction on $i$ and the proof of assertion (1).

For the proof of assertion (2), observe that the process $\sum_{i=1}^k
\delta_{(a^k_i, Y^k_i)}$ converges on $[t,1]\times(c,\infty]$ to
$N=\sum_j \delta_{(\tau_j,Y_j)}$. Further, by induction hypothesis
we have uniform convergence of $\sideset{^{m-1}}{^k}{
{u}}(s,x)$ to $u^{m-1}(s,x)$ as $k\to\infty$. From this, we obtain
convergence of the transformed point processes
\[
\sideset{^m}{^k}{{N}}  = \sum_{i=1}^k \delta_{(a^k_i,
\sideset{^{m-1}}{^k}{{u}}(a^k_i,Y^k_i))} \stackrel{d}{\longrightarrow}  N^m = \sum
_{j} \delta_{(\tau_j,u^{m-1}(\tau_j,Y_j))}, \qquad k\to\infty,
\]
on $M_{u^{m-1}}\cap[t,1]\times\overline\R$ and thus convergence of
the optimal stopping curves of these processes, which proves (2).


Based on (1) and (2), we obtain the estimate
\begin{eqnarray*}
 &&P \bigl(u^m_n(t,x)\ge u^m(t,x) + \epsilon \bigr)
\\
&& \qquad \le P  \biggl(\sideset{^m}{^{n,k}_0}{{W}}(x) \ge \underbrace
{\sideset{^m}{^k}{{u}}(t,x)}_{\sideset{^m}{^k_0}{{u}}
(x)} + \frac{\epsilon}{2}  \biggr) +
P\biggl (u^m(t,x) \le\sideset{^m}{^k}{{u}} (t,x) -\frac
{\epsilon}{2} \biggr).
\end{eqnarray*}
The right-hand side converges for $n\to\infty$ and $k\to\infty$ to
$0$. Thus, we have shown
\[
\lim_{n\to\infty} P\bigl(u^m_n(t,x) \ge u^m(t,x)+\epsilon\bigr) = 0.
\]
To obtain convergence in probability, we next establish that $\liminf
_{n\to\infty} E u^m_n(t,\break x) \ge u^m(t,x)$. This however is
implied by the inequality
\[
Eu^m_n(t,x) \ge E [ X^n_{T_1}\vee\cdots \vee X^n_{T_m} \vee x ]
\]
holding true for all $\mathcal{F}^{n}$-stopping times $tn<T_1<\cdots
<T_m\le n$, and in particular for
\begin{eqnarray*}
\hat{T}^{n,m}_1(t,x) & := & \min\biggl\{ tn < i \le n-m+1 \dvtx  X^n_i > \gamma
^m\biggl(\frac{i}{n},x\biggr)\biggr\},
\\
\hat{T}^{n,m}_\ell(t,x) & := & \min\biggl\{ \hat{T}^{n,m}_{\ell-1}(t,x)
< i \le n-m+\ell\dvtx \\
&&\hphantom{\min\biggl\{}
X^n_i > \gamma^{m-\ell+1}\biggl(\frac{i}{n}, X^n_{\hat{T}^{n,m}_{\ell
-1}(t,x)}\vee x\biggr)\biggr\}
\end{eqnarray*}
for $2\le\ell\le m$. Proposition \ref{prop:5.1} then implies the
above statement.

For $c=-\infty$, we obtain similarly the convergence
$u^m_n(t,x)\stackrel{P}{\longrightarrow} u^m(t,x)$ for $x>-\infty$.
Then the convergence of $u^m_n(t,-\infty) \stackrel
{P}{\longrightarrow} u^m(t)$ results as follows:
\[
u^m_n(t,-\infty)  \le  u_n^m(t,x) \stackrel{P}{\longrightarrow}
u^m(t,x) \downarrow u^m(t) \qquad\mbox{as } x\downarrow-\infty.
\]
This implies that $\lim_{n\to\infty}P(u_n^m(t,-\infty)\ge u^m(t) +
\epsilon) = 0$ for all $\epsilon>0$. Let $\hat{T}^{n,m}_1(t)$ be the
stopping times from condition (L$^m$) and let
\[
\hat{T}^{n,m}_\ell(t) := \min\biggl\{ \hat{T}^{n,m}_{\ell-1} (t) < i
\le n-m+\ell\dvtx  X^n_i > \gamma^{m-\ell+1} \biggl(\frac{i}{n}, X^n_{\hat
{T}^{n,m}_{\ell-1}(t)}\biggr)\biggr\}
\]
for $2\le\ell\le m$. Then we obtain by Proposition \ref{prop:5.1}
and uniform integrability of $(X^n_{\hat{T}^{n,m}_1(t)})_{n\in\N}$ that
\begin{eqnarray*}
 &&E u^m_n(t,-\infty) \ge E\bigl[ X^n_{\hat{T}^{n,m}_1(t)} \vee
\cdots \vee X^n_{\hat{T}^{n,m}_m(t)} \bigr]
\\
&& \qquad \stackrel{n\to\infty}{\longrightarrow} E\bigl[\overline
{Y}_{T^m_1(t,-\infty)} \vee\cdots \vee\overline{Y}_{T^m_m(t,-\infty
)}\bigr] = u^m(t).
\end{eqnarray*}
Thus, $\liminf_{n\to\infty} E u^m_n(t,-\infty)\ge u^m(t)$.
As consequence, we obtain $u^m_n(t,\break -\infty) \stackrel
{P}{\longrightarrow} u^m(t)$ which was to be shown.

(b) For the proof of (b), see \autorF.

(c)
For $c=-\infty$, we obtain the statement using uniform integrability
and Proposition~\ref{prop:5.1}. For $c\in\R$ holds
\begin{eqnarray*}
 &&E  [X^n_{\hat{T}^{n,m}_1} \vee\cdots \vee X^n_{\hat
{T}^{n,m}_m} ]
\\
&& \qquad = E [X^n_{\hat{T}^{n,m}_1} \vee\cdots \vee X^n_{\hat
{T}^{n,m}_m}\vee c ]
\\
&& \qquad  \quad {} - \int_{\{X^n_{\hat{T}^{n,m}_1}
\vee\cdots \vee X^n_{\hat{T}^{n,m}_m}< c\}}  %
 (c - X^n_{\hat{T}^{n,m}_1}\vee\cdots \vee X^n_{\hat{T}^{n,m}_m}
 )\,dP.
\end{eqnarray*}
The first term converges by Proposition \ref{prop:5.1} to the stated
limit. The modulus of the second term can be estimated from above by
\[
\int_{\{X^n_{\hat{T}^{n,m}_m}< c\}}    (c -
X^n_{\hat{T}^{n,m}_m}  )\,dP \le \int_{\{X^n_n< c\}}
   (c-X^n_n ) \,dP \le E|X^n_n-c| \to0.
\]
\upqed
\end{pf}

\begin{rem}
The reason for restricting in (b) to independent sequences is the
necessity to give estimates of $u_n(t,x)$ from above [cf. the case
 $m=1$  in \autorF]. In the dependent case, this amounts to \eqref
{eq:5.4}. For  $m\ge2$  in contrast to the case $m=1$ one has to
consider terms of the form $\max_{s< {i}/{n}\le t}
u_n^{m-1} (\frac{i}{n},X^n_i)$. It seems however difficult to
establish the necessary point process convergence of $\sum_{i=1}^n
\delta_{( {i}/{n},u_n^{m-1}( {i}/{n},X_i))}$ in the general
dependent case.
\end{rem}


\section{Optimal $\mathbf{m}$-stopping of i.i.d. sequences with
discount and observation costs}
\label{sec:6}

As application, we study in this section the optimal $m$-stopping of
i.i.d. sequences with discount and observation costs. In the case
$m=1$, this problem has been considered in various degree of generality
in \citet{KennedyKertz-90}, \citet{KennedyKertz-91}, \autorKRx{-00b} and
\autorFR.

Let $(Z_i)_{i\in\N}$ be an i.i.d. sequence with d.f. $F$ in the
domain of attraction of an extreme value distribution $G$, thus for
some constants $a_n>0$, $b_n\in\R$
%
\begin{equation}\label{eq:6.1a}
n\bigl(1-F(a_nx+b_n)\bigr)\to-\log G(x), \qquad x\in\R.
\end{equation}
Consider $X_i=c_iZ_i+d_i$ the sequence with discount and observation
factors, \mbox{$c_i>0$}, $d_i\in\R$ and both sequences monotonically
nondecreasing or nonincreasing. For convergence of the corresponding
imbedded point processes
%
\begin{equation}\label{eq:6.2a}
\hat N_n=\sum_{i=1}^n \delta_{({i}/{n}, ({X_i-\hat
b_n})/{\hat a_n})}
\end{equation}
the following choices of $\hat a_n$, $\hat b_n$ turn out to be appropriate:
%
\begin{eqnarray}\label{eq:6.3a}
\hat a_n&:=& c_n a_n, \hat b_n:=0  \qquad
\mbox{for } F\in D(\Phi_\alpha) \mbox{ or } F\in D(\Psi_\alpha
),\nonumber
\\[-8pt]
\\[-8pt]
\hat a_n&:=& c_n a_n,  \hat b_n:=c_nb_n+d_n   \qquad \mbox{for }
F\in D(\Lambda),
\nonumber
\end{eqnarray}
where $\Phi_\alpha$, $\Psi_\alpha$, $\Lambda$ are the Fr\'{e}chet,
Weibull, and Gumbel distributions and $a_n$, $b_n$ are the
corresponding normalizations in \eqref{eq:6.1a}. We give further
conditions on $c_i$, $d_i$ to establish point process convergence in
\eqref{eq:6.2a}. Related conditions are given in \citet
{deHaanVerkaade-87} in the treatment of i.i.d. sequences with trends,
respectively, in \autorKR{-00b}\ (\citeyear{KuehneRueschendorf-00b}).

Unlike before, $c$ denotes here a general constant and not the
guarantee value. The guarantee value of $N$ is in case $\Phi_\alpha$
given by 0 and in cases $\Psi_\alpha$, $\Lambda$ given generally by
$-\infty$. This application shows in particular the importance of
treating the case with lower boundary $-\infty$ as in Sections \ref
{sec:2} and \ref{sec:3} of this paper, respectively, in \autorFR. We state the
optimality results for all three cases.
It turns out that in all of the following examples the intensity
functions of the transformed Poisson processes are of
the form studied in Section \ref{sec:4}. Hence, we obtain an explicit
form of the solutions and optimal stopping curves.

We first consider the case of Fr\'{e}chet limits.

\begin{theorem}
\label{theo:6.1}
Let $F\in D(\Phi_\alpha)$ with $\alpha>1$ and $F(0)=0$
(i.e., $Z_i>0$ \mbox{$P$-a.s.}). We assume that $b_n=0$ and
also convergence
\[
\frac{d_n}{c_n a_n} \to d, \qquad \frac{c_{\lfloor tn\rfloor}}{c_n}
\to t^{c}\qquad\forall t\in[0,1]
\]
with constants $c$, $d\in\R$, as well that $c_n$ does not converge to 0.
Assume that $c>-\frac{1}{\alpha}$ and that the function $R\dvtx (d,\infty
)\to\R$,
%
\begin{equation}\label{eq:6.1b}
R(x)  :=  x + \frac{\alpha}{\alpha-1}  \frac{1}{1+c\alpha}
(x-d)^{-\alpha+1}, \qquad x\in(d,\infty),
\end{equation}
has no zero point. Then it holds:
\begin{longlist}[(b)]
\item[(a)]
%
\begin{equation}\label{eq:6.2b}
\frac{E [X_{T_1^{n,m}}\vee\cdots \vee X_{T_m^{n,m}}]}{\hat{a}_n} \to
u^m(0) > 0,
\end{equation}
where $u^m(t)$ is the $m$-stopping curve of the Poisson process $\hat
{N}$ with intensity function
\[
\hat{G}(t,y) = t^{c\alpha} (y-dt^{c+ {1}/{\alpha}})^{-\alpha} =
H \biggl(\frac{y}{v(t)} \biggr) \frac{v'(t)}{v(t)}\qquad\mbox{on }
M_{\hat{f}}.
\]
Here $v(t):=t^{c+ {1}/{\alpha}}$, $H(x):=\frac{\alpha}{\alpha
c+1}(x-d)^{-\alpha}$ and $\hat{f}(t):=dt^{c+{1}/{\alpha}}$.

\item[(b)]
Let $\gamma^1,\ldots ,\gamma^m$ be the functions defined in \eqref
{eq:3.17} for $\hat{N}$. Then
\begin{eqnarray*}
\hat{T}^{n,m}_1   %
&:=&  \min\biggl\{  1\le i \le n-m+1  \dvtx  X_i > \hat
{a}_n\gamma^m \biggl(\frac{i}{n},d\biggr)\biggr\},\\
\hat{T}^{n,m}_\ell  %
&:=&   \min\biggl\{  \hat{T}^{n,m}_{\ell-1} < i \le n-m+\ell
 \dvtx
X_i > \hat{a}_n\gamma^{m-\ell+1}\biggl( \frac{i}{n},  \biggl(\frac
{1}{\hat{a}_n} X_{\hat{T}^{n,m}_{\ell-1}}  \biggr) \vee d\biggr)\biggr\}
\end{eqnarray*}
for $2\le\ell\le m$ are asmptotically optimal sequences of
$m$-stopping times, that is, the limit in \eqref{eq:6.2b} is attained
also for these sequences.
\end{longlist}
\end{theorem}

The next result concerns the Weibull limit case.

\begin{theorem}\label{theo:6.2}
Let $F\in D(\Psi_\alpha)$ with $\alpha>0$ and $F(0)=1$
(i.e., $Z_i\le0$ $P$-a.s.). Further let $a_n\downarrow
0$ and $b_n=0$, and
\[
\frac{d_n}{c_n a_n} \to d, \qquad\frac{c_{\lfloor tn\rfloor}}{c_n}
\to t^{c}  \qquad\forall t\in[0,1]
\]
for constants $c$, $d\in\R$. If $d_n>0$, then assume that either
$(d_n)_{n\in\N}$ is monotonically nondecreasing or $c_n a_n$ does not
converge to 0.
\begin{longlist}[(b)]
\item[(a)]
If $c<\frac{1}{\alpha}$ and $d\le0$, then it holds
%
\begin{equation}\label{eq:6.3b}
\frac{E[ X_{T_1^{n,m}}\vee\cdots \vee X_{T_m^{n,m}}]}{\hat{a}_n} \to
u^m_{c, d}(0) < 0.
\end{equation}

\item[(b)]
If $c>\frac1\alpha$ and the function $R\dvtx \R\to\R$,
%
\begin{equation}\label{eq:6.4a}
R(x) :=
\cases{\displaystyle
 x  ,&\quad if   $x\ge d $,\cr\displaystyle
 x - \frac{\alpha}{\alpha+1}   \frac{1}{1-c\alpha}
(-x+d)^{\alpha+1}  ,&\quad if   $x< d$,
}
\end{equation}
has no zero point then \eqref{eq:6.3b} holds with $u_{c,d}^m(0)>0$.
Here $u_{c,d}^m(t)$ is the $m$-stopping curve of the Poisson process
$\hat N=\hat N_{c,d}$. $\gamma_{c,d}^m$ are the corresponding inverse
functions defined in \eqref{eq:3.17} and \eqref{eq:3.18}.

\item[(c)]
Let $(w_n)$ be an increasing sequence $w_n<0$ such that $n(1-F(w_n))
\to\frac{\alpha+1}{\alpha}$ [e.g., $w_n= - (
\frac{\alpha+1}{\alpha}  )^{{1}/{\alpha}}
a_n$]. Define functions $v_n^m$ by
\[
v_n^m(t) := \frac{\gamma^m_{c,0}(t)}{u_{0,0}(t)}\frac{w_{\lfloor
(1-t)n\rfloor}}{a_n} + \gamma^m_{c,d}(t) - \gamma^m_{c,0}(t),
\]
where $\gamma^m_{c,0}(t)= -\Phi^{m-1} (r_m) u_{c,0}(t)$.
Then the $m$-stopping times defined by
\begin{eqnarray*}
\hat{T}^{n,m}_1 & :=& \min\biggl\{  1\le i \le n-m+1 \dvtx  X_i > \hat{a}_n
v^m_n\biggl(\frac{i}{n}\biggr)\biggr\},
\\
\hat{T}^{n,m}_\ell& :=& \min\biggl\{  \hat{T}^{n,m}_{\ell-1} < i \le
n-m+\ell\dvtx
X_i > \hat{a}_n\gamma^{m-\ell+1}_{c,d}\biggl(\frac{i}{n}, \frac
{1}{\hat{a}_n}X_{\hat{T}^{n,m}_{\ell-1}}\biggr)\biggr\}
\end{eqnarray*}
for $2\le\ell\le m$, are asymptotically optimal, that is, convergence
as in \eqref{eq:6.3b} does also hold for them.
\end{longlist}
\end{theorem}

The final result concerns the Gumbel case.

\begin{theorem}\label{theo:6.3}
Let $F\in D(\Lambda)$ and assume
\[
\frac{b_n}{a_n} \biggl(1-\frac{c_{\lfloor tn\rfloor}}{c_n} \biggr)
\to c\log(t),\qquad\frac{d_n-d_{\lfloor tn\rfloor}}{c_n a_n} \to
d\log(t)\qquad\forall t\in[0,1]
\]
for some constants $c$, $d\in\R$. Assume also that $(c_n)_{n\in\N}$
and $(d_n)_{n\in\N}$ monotonically nondecreasing.
\begin{longlist}[(b)]
\item[(a)]
If $c+d<1$, then
%
\begin{equation}\label{eq:6.4b}
\frac{E[X_{T_1^{n,m}}\vee\cdots \vee X_{T_m^{n,m}}] -\hat{b}_n}{\hat
{a}_n} \to u^m(0),
\end{equation}
where $u^m(t)$ is the $m$-stopping curve of the Poisson process $\hat
{N}$ with intensity function
\[
\hat{G}(t,y) = e^{-y} t^{-(c+d)}\qquad\mbox{on $[0,1]\times\R$}.
\]

\item[(b)]
Let $\gamma^1,\ldots ,\gamma^m$ be the inverse functions defined in
\eqref{eq:3.17} and \eqref{eq:3.18}, let $(w_n)_{n\in\N}$ be an
increasing sequence with $\lim_{n\to\infty} n(1-F(w_n))=1$
(e.g., $w_n:=b_n$). Let $v_n^m$ be defined as
\[
v^m_n(t) := \frac{w_{\lfloor(1-t)n\rfloor}-b_n}{a_n} + \gamma^m(t)
- \log(1-t).
\]
Then
\begin{eqnarray*}
\hat{T}^{n,m}_1 & :=& \min\biggl\{ 1\le i \le n-m+1 \dvtx  X_i > \hat
{a}_nv_n^m\biggl(\frac{i}{n}\biggr)+\hat{b}_n\biggr\},
\\
\hat{T}^{n,m}_\ell& :=& \min\biggl \{ \hat{T}^{n,m}_{\ell-1} < i
\le n-m+\ell\dvtx
\\
&&\hphantom{\min\biggl \{}   X_i > \hat{a}_n\gamma^{m-\ell+1}
\biggl(\frac{i}{n}, \frac{X_{\hat{T}_{\ell-1}^{n,m}} -\hat{b}_n}{\hat
{a}_n}  \biggr) +\hat{b}_n
 \biggr\}
\end{eqnarray*}
define an asymptotic optimal sequence of $m$-stopping times, that is,
convergence as in \eqref{eq:6.4b} holds for them.
\end{longlist}
\end{theorem}

For details of the proof, we refer readers to \autorFx{pages 75--77}.

%

\printaddresses

\end{document}